\magnification=1100

\font\titfont=cmr10 at 16 pt

\font\headfont=cmr10 at 12 pt     
 
 \font\aufont=cmr10 at 14pt

\font \fr = eufm10

\overfullrule=0in

\def\boxit#1{\hbox{\vrule
 \vtop{%
  \vbox{\hrule\kern 2pt %
     \hbox{\kern 2pt #1\kern 2pt}}%
   \kern 2pt \hrule }%
  \vrule}}

  \def\harr#1#2{\ \smash{\mathop{\hbox to .3in{\rightarrowfill}}\limits^{\scriptstyle#1}_{\scriptstyle#2}}\ }

\def\wee{\wedge\cdots\wedge}

\def\ss{\subset}
\def\sse{\subseteq}
\def\half{\hbox{${1\over 2}$}}
\def\smfrac#1#2{\hbox{${#1\over #2}$}}
\def\oa#1{\overrightarrow #1}

\def\log{{\rm log}}
\def\Hess{{\rm Hess}}

\def\tr{{\rm tr}}

\def\det{{\rm det}}

\def\arr{\longrightarrow}
\def\supp{{\rm supp}}

\def\rn{\bbr^n}

 \def\cd{{\cal C}}
 
 \def\cm{{\cal M}}

\def\Theorem#1{\medskip\noindent {\bf THEOREM \bf #1.}}
\def\Prop#1{\medskip\noindent {\bf Proposition #1.}}
\def\Cor#1{\medskip\noindent {\bf Corollary #1.}}
\def\Lemma#1{\medskip\noindent {\bf Lemma #1.}}
\def\Remark#1{\medskip\noindent {\bf Remark #1.}}
\def\Note#1{\medskip\noindent {\bf Note #1.}}
\def\Def#1{\medskip\noindent {\bf Definition #1.}}

\def\Ex#1{\medskip\noindent {\bf Example \bf    #1.}}

\def\pf{\medskip\noindent {\bf Proof.}\ }
\def\qed{\hfill  $\vrule width5pt height5pt depth0pt$}
\def\mathqed{\vrule width5pt height5pt depth0pt}

\def\qedqed{\hfill  $\vrule width5pt height5pt depth0pt$ $\vrule width5pt height5pt depth0pt$}

\def\hk{\_{\rm l}\,}

\def\w{\wedge}

   \def\cc{{\cal C}}     
      
   \def\co{{\cal O}}
\def\ce{{\cal E}}   
\def\ch{{\cal H}}   \def\cm{{\cal M}}
\def\cs{{\cal S}}   
\def\cd{{\cal D}}
\def\cl{{\cal L}}

\def\cf{{\cal F}}

\def\gerM{{\fr{\hbox{M}}}}

\def\vf{\varphi}

\def\wt{\widetilde}

\def\and{\qquad {\rm and} \qquad}
\def\arr{\longrightarrow}
\def\ol{\overline}
\def\bbr{{\bf R}}\def\bbo{{\bf O}}
\def\bbc{{\bf C}}

\def\bbz{{\bf Z}}
\def\bbp{{\bf P}}
\def\bbd{{\bf D}}

\def\a{\alpha}
\def\b{\beta}

\def\e{\epsilon}

\def\g{\gamma}

\def\l{\lambda}
\def\o{\omega}

\def\s{\sigma}
\def\x{\xi}
\def\z{\zeta}

\def\D{\Delta}
\def\L{\Lambda}
\def\G{\Gamma}
\def\O{\Omega}

\def\dbar{\ol{\partial}}

\def\bo{\partial \Omega}

\def\fa{{\rm\ \  for\ all\ }}

\def\M{{\bf M}}

 \def\LAG{{\rm LAG}}

\def\bbt{{\bf T}}
 \def\bre{\bar e}
 \def\dn{{\bbd^n}}
 \def\D{4}
 \def\SS{5}
 \def\DD{6}
 \def\R{7}
 \def\E{8}
 \def\F{9}
 \def\G{10}
 \def\H{11}
 \def\QQ{13}
 \def\MM{14}
 \def\UU{16}
 \def\TT{15}
 \def\LAG{{\rm LAG}}
  \def\LAGt{\wt{\rm LAG}}
  \def\SLAGt{\wt{\rm SLAG}}
 \def\bra#1#2{{\langle #1, #2\rangle}}
 \def\Sph{{\bf S} \!\!\! \! /}
 \def\be{{\bf e}}
 \def\Ta{{\overrightarrow T}}
 \def\Gs{{G_{\rm space}(n,\dn)}}
 \def\Gsp{{G_{\rm space}^+(n,\dn)}}
 \def\Gsn{{G_{\rm space}^-(n,\dn)}}
 \def\CI{{\cal I}}

\centerline{\titfont  Split Special Lagrangian Geometry}
\medskip
 \vskip .4in

\centerline{\aufont F. Reese Harvey and H. Blaine Lawson, Jr.$^*$}
\vglue .9cm
\smallbreak\footnote{}{ $ {} \sp{ *}{\rm Partially}$  supported by 
the N.S.F. }

\vskip .5in
\centerline{\bf ABSTRACT} \bigskip
  \font\abstractfont=cmr10 at 10 pt
{{\parindent= .53in
\narrower\abstractfont \noindent

One purpose of this article is to draw attention to the seminal work of J. Mealy in 1989 on calibrations
in semi-riemannian geometry where split SLAG geometry was first introduced.  The natural 
setting is provided by doing geometry with the  complex numbers $\bbc$ replaced by  the double numbers $\bbd$,
where $i$  with $i^2=-1$ is replaced by $\tau$ with $\tau^2=1$.  
A rather surprising amount of complex geometry carries over, almost untouched,
and this has been the subject of many papers.  We briefly review this material and, in particular,
we discuss Hermitian $\bbd$-manifolds with trivial canonical
bundle, which provide the background space for the geometry of split SLAG submanifolds.

 A removable singularities result is proved for split SLAG subvarieties.  It implies,
in particular, that there exist no split SLAG cones, smooth outside the origin, other than 
planes.  This is in sharp contrast to the complex case.

Parallel to the complex case, space-like Lagrangian submanifolds are stationary
if and only if they are $\theta$-split SLAG for some phase angle $\theta$, and infinitesimal deformations
of split SLAG submanifolds are characterized by harmonic 1-forms on the submanifold.

We also briefly review the recent work of Kim, McCann and Warren who have shown that split 
Special Lagrangian geometry is directly related to the Monge-Kantorovich mass transport
problem.

}}

\vfill\eject

\centerline{\bf TABLE OF CONTENTS} \bigskip

{{\parindent= .1in\narrower\abstractfont \noindent

\qquad 1. Introduction.\smallskip

\qquad 2.      Double Numbers and Double Holomorphic Functions.
\smallskip

\qquad 3.   Double $n$-Space. \smallskip

\qquad 4.   The Special Lagrangian Calibration.
\smallskip

\qquad 5.     Split SLAG Submanifolds.
  
\smallskip

\qquad 6.    Split SLAG Graphs in D$^n$
  
\smallskip

\qquad 7. Removable Singularities for   split  SLAG  Varieties.

\smallskip

\qquad 8. Double Manifolds.
   
\smallskip

\qquad 9.  Hermitian D-Manifolds.
   
\smallskip

\qquad 10.  K\"ahler D-Manifolds.
   
\smallskip

\qquad 11.  Ricci-Flat  K\"ahler D-Manifolds.

\smallskip

\qquad 12. Split SLAG Submanifolds in the General Setting.

\smallskip

\qquad 13.  Deformations and the McLean Theorem in the Split Case

\smallskip

\qquad 14.  Relation to the Mass Transport Problem -- 
Work of Kim, McCann and Warren

\smallskip

\qquad 15. Further Work of Mealy.

\smallskip

\qquad 16.    Lagrangian Submanifolds of Constant Phase and Volume Maximization

}}

\vskip .3in

{{\parindent= .3in\narrower

\noindent
{\bf Appendices: }\medskip

A. \  A Canonical Form for Space-like $n$-Planes in D$^n$.
\smallskip

B. \  Projective and Hermitian Projective Varieties over D.
\smallskip

C. \ Degenerate Projections.
\smallskip

D. \ Singularities and Semi-Riemannian Calibrations.
\smallskip

E. \ Singularities in Split SLAG Geometry.
\smallskip

F. \ Split SLAG Singularities in Dimension Two.
\smallskip

}}

\vfill\eject


   \vskip .3in
 

\centerline{\headfont \ 1.\  Introduction. }  
\medskip

In an interesting recent series of papers [KM], [W$_1$], [KMW]  Y.-H. Kim, R. McCann and M. Warren
established a direct relationship between the classical Monge-Kantorovich mass transport problem
and split special Lagrangian geometry. This geometry is one of the basic calibrated semi-riemannian
geometries introduced by Jack Mealy in his seminal work [M$_{1,2}$] in 1989-90.
The ambient spaces are analogues of complex K\"ahler manifolds for which the literature is
vast and not unified$^1$
\footnote{}{ $ {} \sp{1}{\rm For \ }$ {example, see [AMT], [CFG], [EST], [GM], [GGM], [IZ], [SS] and  references cited therein.}}
 -- they are sometimes called paraK\"ahler and sometimes  Bi-Lagrangian (we shall call
them K\'ahler $\bbd$-manifolds). 
 Despite the widespread interest, 
Mealy's  results about calibrations on such manifolds seem to have gone unnoticed. 
One  purpose of this article is to
give a brief self-contained discussion of $\bbd$-manifolds
so that we can discuss split special Lagrangian geometry and present the
work of Mealy in this important case.
  We shall then explain how the work of Kim-McCann-Warren 
fits into the picture. Along the way we will also establish some new results concerning split special Lagrangian
subvarieties. One is a removable-singularities result which implies, in particular, 
that any split SLAG cone, which is smooth and connected outside the origin, is a plane
(in sharp contrast to the classical non-split case where such cones exist in abundance).
 Another can be interpreted as the strongest possible  regularity result for  $d$-closed split special
Lagrangian rectifiable currents in dimension 2.

Throughout the exposition we  emphasize the two distinct
 approaches to $\bbd$-manifolds.  The first is in strong analogy with complex geometry,
where the double numbers $\bbd$ (see \S 2) are considered as a replacement for  the complex numbers
$\bbc$, and a surprising    number of definitions and results carry over. 
In the second viewpoint the double numbers are viewed as the algebra $\bbd=\bbr\oplus\bbr$.
Here one focuses on the canonically
defined pair of transverse $n$-dimensional foliations which give rise locally to null coordinates.
In this setting many classical objects, such as  holomorphic $n$-forms, which are written as usual in complex
notation, appear quite differently in null coordinates.
A striking example is the equation for the (split) special Lagrangian potential function.
In ``complex coordinates''  the equation resembles that of the non-split  case -- all the odd
elementary symmetric functions of the hessian must sum to zero.  In null coordinates, it  becomes  the 
Monge-Amp\`ere equation $\det\, \Hess\, u = 1$.

The algebra $\bbd$ mentioned in the paragraph above has at least four names: {\sl the double numbers,
the para-complex numbers, the Lorentz numbers, and the hyperbolic numbers.}
It is the  algebra over $\bbr$ generated by 
1 and $\tau$ where $\tau^2=1$. In other words $\bbd$ is the Clifford algebra for the 
negative definite quadratic form on $\bbr$ (while $\bbc$ comes from the positive definite one).
 The algebras $\bbc$ and $\bbd$ are the only commutative normed algebras aside from $\bbr$. They are like 
twin sisters living in parallel universes, one elliptic and the other hyperbolic. Writing elements
in $\bbd$ as $z=x+\tau y$, and setting $\bar z=x-\tau y$ one can proceed formally to define
$\bbd$ holomorphic functions (of one and several variables), and the analogues of complex manifolds
(called $\bbd$-manifolds here),
holomorphic vector bundles, the Dolbeault complex, and much more.  There are $\bbd$ projective spaces
and grassmannians. In fact every real algebraic variety has an associated $\bbd$-variety via base change.
Even more striking are the metric analogues. There are hermitian $\bbd$-manifolds, K\"ahler  $\bbd$-manifolds
and  Calabi-Yau $\bbd$-manifolds.  Many basic theorems carry over.  This is discussed below.

However, as the reader may have noticed, after setting $e=\half(1+\tau)$ and $\bar e=\half(1-\tau)$,   $\bbd$ 
becomes the double algebra $\bbr\oplus\bbr$ of elements
$u e + v \bar e$ where $u,v\in\bbr$ (and $e^2=e, {\bar e}^2=\bar e$ and $e\bar e=0$).  
These are the null-coordinates mentioned above. The analogues of  many fundamental objects look quite different in these coordinates.
For example, $-\tau \partial\dbar = d_u d_v$ where $d_u$ and $d_v$ are the standard deRham differentials
in the variables $u$ and $v$.

The paper is organized as follows.  After a brief introduction to analysis over $\bbd$ we examine the
split special Lagrangian calibration $\Phi={\rm Re}\,dz$ in $\bbd^n$ and prove Mealy's  sharp inequality which 
states that $\Phi$ is a (reverse) calibration which is minimized on space-like special Lagrangian $n$-planes.
By the fundamental lemma of calibrations this implies that 
split SLAG submanifolds are homologically volume-maximizing. This is done in \S\S 1-5.

Recall that the natural exterior differential system (EDS) for SLAG submanifolds of $\bbc^n$
is generated by Im$\,dz$ and the standard symplectic (K\"ahler) form  $\o$.  A submanifold
is an integral submanifold for this EDS if and only if it is  SLAG when correctly oriented.
The parallel to this in the split case (Theorem 5.2) is almost identical.  Again the differential system
is generated by Im$\,dz$ and $\o$ (this time on $\dn$), but an integral submanifold is split SLAG
modulo orientation if and only if it is space-like.  These more general integral submanifolds, which are
not necessarily space-like, will be called {\sl unconstrained} split SLAG submanifolds

In \S 6 we examine split SLAG   submanifolds $M$ which are graphs of mappings.
If $M$ is the oriented graph of $F:\O\to {\rm Im}\, \dn$ where $\O$ is a simply-connected domain
in ${\rm Re} \,\dn \cong \rn$, then  following Mealy, $F= df$ for a potential function $f$ for which 
$$
\quad
{\rm Im} \, \det_\bbd (I+\tau \Hess\, f)\ =\ 0 \qquad \quad {\rm and \ \qquad} -I \ < \  \Hess f \ <\  I.
\eqno{(1.1)}
$$
If on the other hand, $M$ similarly is expressed as the graph $v=G(u)$ over a domain in null coordinates,
then $G=dg$ for a potential function $g$ satisfying (cf. [Hi$_1$])
$$
\quad
\det\,\Hess\, g \ =\ 1  \quad  \qquad {\rm and \ \qquad } \Hess \, g \ >\ 0.
\eqno{(1.2)}
$$

In \S 7 we prove the following removable-singularities result.

\Theorem{\R.1} {\sl
Let $\O\ss\rn$ be a convex domain in $\rn$ and consider the 
``tube'' domain $\Omega\times \tau\rn \ss  \rn\oplus \tau \rn=\dn$.
Let $\Sigma\ss\O$ be a compact subset of Hausdorff $(n-2)$-measure zero. Then any closed, 
connected   split  SLAG  submanifold 
$$
M\ss (\O-\Sigma)\times \tau\rn
$$ 
has closure in $\O\times \tau\rn$ which is the graph of a real analytic mapping
$F:\O\to\rn$.  Furthermore, $F= d f$ where $f:\O\to \bbr$ satisfies 
(1.1).
}

 \Cor{\R.2. (Absence of Cones)} {\sl 
 Suppose $C\ss \bbd^n$  is a   split SLAG cone which is regular and connected  outside the origin.
 Then $C$ is an $n$-plane.
 }\medskip

  This is in sharp contrast to the non-split case where ${\rm SLAG}$-cones exist
in all dimensions and there are many ${\rm SLAG}$-varieties with singularities of high codimension
(cf. [HL$_1$],  [J$_*$],[H$_*$], [HK$_*$]  and the references therein).

The proof of Theorem \R.1 relies on playing off the two different coordinate pictures for $z\in \dn$.
With $z=x+\tau y$
we show that the projection of $M$ to $\O-\Sigma$ must be
a covering map. This leads to a potential with a Lipschitz extension across $\Sigma$.
 We then rotate to null  coordinates  $z=ue+v\bar e$ (the Cayley transform)
and apply deep results of Caffarelli [Ca$_*$] to get full regularity.

Corollary \R.2 would suggest that singularities of codimension $\geq 3$ do not exist
on split SLAG subvarieties. By contrast, singularities of codimension-2 certainly do exist
in all dimensions.  It turns out that when $n=2$ the exterior differential system for (unconstrained) 
split SLAG submanifolds is equivalent to the Cauchy-Riemann system for complex curves.
Specifically we have the following.

\vfill\eject

\Theorem{F.3} {\sl
There exists a permutation of real coordinates $\bbc^2 \leftrightarrow \bbd^2$, such that:
\medskip
\hskip .2in 
(a)\  \ holomorphic chains
are transformed to unconstrained split SLAG currents

\hskip .5in  and vice-versa.
\medskip

\hskip .2in 
(b)  \  \ holomorphic chains  satisfing a 45$^o$-rule 
are transformed to  split SLAG currents

\hskip .5in   and vice-versa.
}\medskip

Since complex curves have local uniformizing parameters, this can be interpreted
as  a  very strong regularity result  for split SLAG subvarieties when $n=2$.
(A split SLAG subvariety is a $d$-closed rectifiable current $T$
whose tangent planes $\Ta$ are space-like special Lagrangian almost everywhere.  
See Appendices D and E for a fuller discussion.)
Theorem F.3  follows from  work of the first author and B. Shiffman [HS], [S]  which 
states that: a $d$-closed rectifiable 2-current defined in an open subset $X$ of $\bbc^2$, whose
  (unoriented) tangent space $\Ta$   is a complex line $\|T\|$-a.e.,  is a  $T$  holomorphic chain
 in $X$. That is, $T$ is a locally finite sum $T=\sum_j n_j [V_j]$ where for each $j$, $n_j\in \bbz$ and 
$V_j$ is a 1-dimensional complex analytic subvariety of $X$.

 This result also clearly  implies 
that isolated singularities {\bf do occur} on split SLAG subvarieties when $n=2$.
They are the singular points of the underlying holomorphic curves.
Taking products of these with appropriately chosen linear subspaces gives
codimension-2 singularities on split SLAG subvarieties in all dimensions $n$.

It is interesting to note that since locally irreducible holomorphic curves have well defined tangent 
lines at singular points, this procedure does not yield interesting cones.  So far, no non-planar
irreducible cones are known.

In sections \E\ through \H\ we discuss the geometry of $\bbd$-manifolds. In particular we 
examine hermitian $\bbd$-manifolds, K\"ahler $\bbd$-manifolds, and Ricci-Flat K\"ahler
$\bbd$-manifolds.  A number of examples are given.
In section  12 we discuss    split SLAG submanifolds in this general setting.
In section 13 we consider the deformation theory for these submanifolds in
the Ricci-flat K\"ahler case.  It is shown that the basic results of McLean [Mc$_{1,2}$]
for the classical (non-split) special Lagrangian submanifolds carry over to this case
(cf. Warren [W$_2$]).

Then in \S 14 we discuss the work of Kim-McCann-Warren on the relevance
of split special Lagrangian geometry 
to the mass transport problem.

In section \TT \ we briefly survey the work of Mealy on other calibrations in geometries of
indefinite signature.  

In section \UU\ we analyze the case of Lagrangian submanifolds of $\dn$ (or more
generally of Ricci-flat K\"ahler $\bbd$-manifolds) on which the restricted metric 
is everywhere non-degenerate of signature $p,q$ ($p+q=n$). We show that
any such manifold $M$ is stationary  (mean curvature zero)  if and only if the restriction of $dz$ to $M$ is of
constant phase (see  [HL$_1$], [D] for parallel cases). 
When $p=n$ (respectively, 0), a compact  oriented constant-phase submanifold
maximizes volume among all oriented space-like (respectively, time-like) submanifolds
with the same boundary.  For $0<p<n$,  volume is not maximized in general,  
but  it is  maximized among oriented Lagrangian submanifolds
of the same type.  In all cases the volume inequality becomes an equality   only
when the competitor is also of the same constant phase.

Appendix A gives a canonical form for a  space-like $n$-planes in $\dn$,
under the action of the $\bbd$-unitary group, which is 
used in proving Mealy's theorem.

Appendix B treats projective $\bbd$-subvarieties and their induced  split 
K\"ahler geometry.

Appendix C shows that the result proved in [HL$_1$] characterizing 
special
Lagrangian submanifolds with  degenerate projections, carries over
to the split case.

Appendix D presents the general theory of calibrations on semi-riemannian manifolds.  
This entails, among other things, the establishing
of a canonical polar form for space-like currents which are representable by integration.
The Fundamental Theorem for Semi-Riemannian
Calibrations is then established for general $\phi$-subvarieties (rectifiable $\phi$-currents).

Appendix E examines the special case of split SLAG subvarieties which are shown to be
 volume maximizing (cf. [KMW]).

Appendix F is devoted to the case $n=2$ which was discussed above.

We would like to thank Luis Caffarelli, Young-Heon Kim and Micah Warren for helpful comments.

   \vskip .3in

\centerline{\headfont \ 2.\  Double numbers and Double Holomorphic Functions. }  
\medskip
In this section we recall the basics of {\sl double} or {para-complex} analysis.
By the {\sl double numbers} we  mean the 2-dimensional algebra  $\bbd$ over $\bbr$
generated by $1$ and $\tau$ with $\tau^2=1$. 
This is the only commutative normed algebra other than the real and complex numbers.
In analogy with the complex numbers,
each $z\in \bbd$ can be written as 
$$
z\ =\ x+\tau y
\eqno{(2.1)}
$$
with $x,y\in \bbr$ defined to be the {\sl real} and {\sl  imaginary parts}.  However, choosing the  basis $e=\half (1-\tau)$ and $\bar e = \half (1+\tau)$
each $z\in \bbd$ can be written as 
$$
z\ =\  ue+ v\bar e
\eqno{(2.2)}
$$
where $(u,v)$ are called the {\sl null-coordinates} of $z$. Note that $e^2=e$, ${\bar e}^2=\bar e$
and $e\bar e=0$. Thus we see that $\bbd$ is just the algebra of diagonal $2\times 2$-matrices
$(u,v)$, i.e., $\bbd = \bbr\oplus\bbr$ as algebras.
Note the relations
$$
\tau\cdot e\ =\ -e\and \tau\cdot \bar e\ =\  \bar e.
\eqno{(2.3)}
$$

The representation (2.1) leads to a vast development of analysis and geometry
in parallel with the complex case. This development is frequently demystified by using
the null coordinates (2.2).

Conjugation in $\bbd$ is defined exactly as in the complex case:  $\bar z=x-\tau y = ve+u\bar e$.
However 
$$
\langle z,z\rangle \ =\ z\bar z = x^2-y^2 = uv
$$
defines a  quadratic form of signature  $(1,1)$.  The algebra $\bbd$ is normed in the sense that
$\langle zw,zw\rangle = \langle z,z\rangle\langle w,w\rangle$.
The notion of vanishing or being zero in complex analysis 
is replaced by the notion of being {\sl null}, i.e., $z\bar z=0$, in which case having an inverse is
impossible. Otherwise, $z^{-1}=  \bar z/\langle z,z\rangle$ as in the complex case.

The multiplicative group $\bbd^*$ of non-null numbers has four connected components.
Let $\bbd^+$ denote the component containing 1.  Exponentiation gives an isomorphism:

\noindent
${\rm exp}:\bbd \harr{\cong}{} \bbd^+$  with inverse $\log: \bbd^+ \to \bbd$. Note that using (2.1) and (2.2)
$$
\eqalign
{\exp(z) \ & =\ \exp(x)\left( \cosh y + \tau  \sinh y\right)  \cr
&=\ e \exp(u) + \bar e \exp(v)  \cr
}
$$
The {\sl  space-like unit sphere} has two components parameterized by $\pm \exp(\tau\theta)$.

As in complex analysis we now treat $\bbd$ as $\bbr^2$ and consider  smooth $\bbd$-valued
functions on open sets $U\ss\bbd$.  This leads to several parallels.
For example the $\bbd$-valued 1-forms 
$$
dz\ =\ dx+\tau dy\ =\ e du+\bar e dv
\and
d\bar z\ =\ dx-\tau dy\ =\ e dv+\bar e du
$$
have duals
$$
 {\partial\over \partial z} \ =\ {1\over 2} \left ( {\partial\over \partial x} +\tau {\partial\over \partial y} \right ) \ =\ 
e {\partial\over \partial u} +\bar e {\partial\over \partial v}
\and
 {\partial\over \partial \bar z}\ =\ {1\over 2} \left ( {\partial\over \partial x} -\tau {\partial\over \partial y} \right ) \ =\ 
e {\partial\over \partial v} +\bar e {\partial\over \partial u}.
$$
However, 
$$
4{\partial^2 \over \partial z\partial\bar z}\ =\ 
{\partial^2 \over \partial x^2} - {\partial^2 \over \partial y^2}\ =\ 
4{\partial^2 \over \partial u\partial v}\
\eqno{(2.4)}
$$
is the wave equation, so that the regularity obtained in the complex case is lacking.
Nevertheless,  we define a smooth function $F:U\to \bbd$ to be {\sl $\bbd$-holomorphic} if 
$$
{\partial F\over \partial \bar z} \ =\ 0.
\eqno{(2.5)}
$$
If we write $F=ef+\bre g$ with $f$ and $g$ real-valued, then using null coordinates,
$$
\dbar  F\ \equiv\ {\partial F\over \partial \bar z} d\bar z\quad{\rm becomes\ \ }
\dbar (ef+\bar eg) \ =\ e\,d_vf + \bar e\, d_u g. \ \ \ {\rm Also,}
$$
$$
\partial  F\ \equiv\ {\partial F\over \partial  z} dz\quad{\rm becomes\ \ }
\partial (ef+\bar eg) \ =\ e \,d_uf + \bar e \,d_v g.
$$
Here    $d_u$ and $d_v$ are the usual exterior differentiation
 operators in the  null coordinates $u$ and $v$.  Thus
\smallskip
\centerline{$F\ =\ ef+ \bar e g$\ is  $\bbd$-holomorphic \qquad$\iff$ } 

\centerline{
$ f=f(u)$ and $g=g(v)$ are functions of the single variable $u$ and $v$ respectively. \hfil (2.6)
 }
 \bigskip\noindent
 For any $F = ef+\bar e g$ we have $2 \, {\rm Re}  F = f+g$ and $2\, {\rm Im}  F= f-g$. Consequently, 
 we see that  the real parts 
  of $\bbd$-holomorphic functions are functions of the form $f(u)+g(v)$.
Observe that 
 just as in the complex  case:
$$
F \ {\rm is\ } \bbd  \ {\rm holomorphic\ and \ Re}\,F =0 \quad\Rightarrow\quad
 F = \tau c \ {\rm is \ constant} 
 \eqno{(2.7)}
 $$
 despite the hyperbolic character of (2.4).
 Finally note   that $$-\tau\partial\dbar = d_u d_v$$ is a real operator. 
 One can show that any function in the kernel of this operator is locally 
 the real part of a $\bbd$-holomorphic function.

   \vskip .3in
 

\centerline{\headfont \ 3.\  Double n-Space D$^n$.}  
\medskip

The standard formulas in several complex variables have rather obvious analogues.
Define
$$
\partial  \ =\ \sum_{k=1}^n dz_k\wedge {\partial \over \partial z_k}
\and
\dbar   \ =\ \sum_{k=1}^n d\bar z_k \wedge {\partial \over \partial \bar z_k}
\eqno{(3.1)}
$$
operating on $\bbd$-valued differential forms (forms with coefficients in $\bbd$ at each point).
These operators satisfy the usual relations: $\partial^2=\dbar^2 = \partial\dbar+\dbar\partial =0$
and $d=\partial + \dbar$. We can define the real operator
$
d^c = \tau(\partial-\dbar) = 2{\rm Im}(\partial)
$
so that $\partial = \half(d+\tau d^c)$. Then $-{\tau\over 2} \partial\dbar ={1\over 4} d d^c$ is a real operator.
So far these are exactly the same (no sign changes) as in the complex case.
Switching to null coordinates the operators are given by
$$
\partial \ =\ e\,d_u + \bar e\, d_v \and \dbar \ =\ e\,d_v + \bar e\, d_u.
\eqno{(3.2)}
$$
 The standard symplectic form $\o$ on $\bbr^{2n}$ can be written in the two ways
 $\o = \sum_j dx_j\wedge dy_j = \half \sum_j du_j\wedge dv_j$
 and has the $-\half \tau \partial\dbar = \half d_u d_v$ potential $\vf = \sum_j z_j{\bar z}_j = \sum_j u_j \bar v_j$.

 \Def{3.1} A $\bbd$-valued function $F$ on an open subset of $\bbd^n$ is $\bbd$-{\sl holomorphic}
 if $\dbar F=0$. Note that in null coordinates $(u,v)$ if we write  $F=e\,f+\bar e\, g$ with $f$ and $g$ real-valued,
  then

 \medskip
\centerline{$F$ is $\bbd$-holomorphic \ $\iff \ d_vf=0$ and $d_ug=0$ \ $\iff\  f=f(u)$ and $g=g(v)$.}
\medskip
 
 On $\bbd^n$ we have the standard $\bbd$-valued hermitian form and real inner product:
 $$
 (z,w) \ =\ \sum_{j=1}^n z_j \bar w_j
 \and
 \langle z,w \rangle \ =\  {\rm Re}(z,w) \ =\  {\rm Re}\sum_{j=1}^n z_j \bar w_j.
 $$
The associated real quadratic form  $\langle z,z \rangle = \sum (x_j^2-y_j^2)$ has split or neutral signature.
 Note that the standard  $\bbd$-valued  hermitian form has imaginary part $-\o$. Hence, exactly as in the complex case 
 $$
 (\cdot,\cdot) \ =\ \langle \cdot,\cdot\rangle - \tau\o (\cdot,\cdot)
 \eqno{(3.3)}
 $$
 The analogue of the complex structure operator $J$ is the operator $Tz=\tau z$ of multiplication
 by $\tau$ (see (8.3)).  Now (3.3) leads to an array of formulas such as:
 $\bra  {z'} {Tz} =\o(z',z)$ relating 
 $\langle \cdot,\cdot\rangle$, $\o$ and $T$  (see (9.3)).

The $\bbd$-linear maps from $\bbd^n$ to $\bbd^m$  correspond to the set $M_{n,m}(\bbd)$ of  $n\times m$-matrices
with entries in $\bbd$.  By standard algebra each square matrix  $A\in M_{n}(\bbd)$ has a $\bbd$-determinant with the 
usual properties such as:
$
\det_\bbd (AB) = (\det_\bbd A)(\det_\bbd B)$,
$
\det_\bbd A^t = \det_\bbd A$, 
 and 
 $
  A\widetilde A^t\ =\ \left(\det_\bbd A\right) I
$
where $\widetilde A$ is the cofactor matrix. 
Thus $A$ has an inverse $A^{-1} \in M_n(\bbd)$ if and only if $\det_\bbd A$ has in inverse, i.e.,  
$\det_\bbd A \in \bbd^*$  is non-null.   
Let ${\rm GL}_n(\bbd) \ss M_n(\bbd)$ denote the group of invertible elements.
Setting $dz \equiv dz_1\wedge\cdots\wedge dz_n$
we have that if $z'=Az$ for  $A\in M_{n}(\bbd)$, then $dz'= \det_\bbd A \,  dz$.
Since $dz\wedge d\bar z$ is a non-null  real multiple of $dx\wedge dy$, this proves that (as in the complex case)
$$
\left(  \det_\bbd A \right)\left( \overline{ \det_\bbd A }\right)\ =\ \det_\bbr A.
$$
Writing $A= eB+\bar e C\in M_{n}(\bbd)$ with $B,C \in M_{n}(\bbr)$ and noting that
$
dz\ =\ e du + \bar e dv
$
yields the simple formula
$$
\det_\bbd(eB+\bar e C)\ =\ e\,\det_\bbr B+\bar  e\,\det_\bbr C
$$
for the $\bbd$-determinant.  In particular, $\det_\bbd A\in \bbd^+ \ \iff\  \det_\bbr B>0 \ {\rm and} \  \det_\bbr C>0$.

Each real linear map from $\bbd^n$ to $\bbd^m$ decomposes into the sum of a 
$\bbd$-linear map and an anti-$\bbd$-linear map.   Thus the Jacobian of a smooth mapping $F$
from (an open subset of)   $\bbd^n$ to $\bbd^m$ decomposes as
$$
J_\bbr(F)\ =\ J_\bbd^{1,0}(F) + J_\bbd^{0,1}(F)\ \equiv \ {\partial F\over \partial z} +  {\partial F\over \partial \bar z}
$$
and $F$ is said to be {\sl $\bbd$-holomorphic} if the anti-linear part $J_\bbd^{0,1}(F)=0$,
i.e., if $J_\bbr(F)$ is $\bbd$-linear at each point. 

In order for a function $F$ between open subsets of $\bbd^n$ to be 
bi-holomorphic we must require that $F$ be holomophic and that 
$\det_\bbd {\partial F\over \partial z}\in \bbd^*$ be non-null, so that 
${\partial F\over \partial z}$ is invertible.

The $\bbd$-unitary group, denoted $U_n(\bbd)$, can be defined, as in the complex case,
in several equivalent ways.  Given $A\in M_n(\bbd)$, let $A^* =  {\overline A}^t$, denote the
conjugate transpose.   Let $e_1,...,e_n, Te_1,...,Te_n$ denote the standard basis for $\bbd^n 
= \rn\oplus\tau\rn$.  A set of vectors $\e_1,...,\e_n$  in $\bbd^n$ is a {\sl space-like $\bbd$-unitary basis
for} $\dn$ if $\e_1,...,\e_n,T\e_1,...,T\e_n$ is a real orthonormal basis with $\langle \e_j, \e_j\rangle =1$
$\langle T\e_j,   T \e_j\rangle = -1$ for all $j$.  A matrix  $A\in M_n(\bbd)$ is called {\sl  $\bbd$-unitary}
if one of the following equivalent conditions holds:
\smallskip

(1)\ \ $(Az, Az)  = (z,z) \ \ $ for all $z\in\dn$

\smallskip

(2)\ \ $A A^* =I$  \ or\ $A^* A =I$ \hfill (3.4)

\smallskip

(3)\ \ $Ae_1,...,Ae_n$ is a space-like $\bbd$-unitary basis for $\dn$.

\smallskip
 
When $n=1$, the unitary group ${\rm U}_1(\bbd)$ is the space-like sphere
$\{\pm e^{\tau\theta} : \theta \in \bbr\}$ with two components.  In general,
${\rm U}_n(\bbd)$ has two components determined by $\det_\bbd A= \pm e^{\tau\theta}$.  
Let ${\rm U}_n^+(\bbd)$ denote the identity component.  Computing in null coordinates, first express
$A\in M_n(\bbd)$ as $A= e \,B + \bar e \,C$ with $B,C \in M_n(\bbr)$, then we see that:
$A\in {\rm U}_n(\bbd)   \iff A= e\, B + \bar e\, (B^t)^{-1}$ for some $B\in {\rm GL}_n(\bbr)$, and that 
$A\in {\rm U}_n^+(\bbd)  \iff  B\in {\rm GL}^+_n(\bbr)$.  Thus ${\rm U}_n^+(\bbd) \cong {\rm GL}^+_n(\bbr)$.

The {\sl special unitary group}  ${\rm SU}_n(\bbd)$ is defined by the additional condition: $\det_\bbd A=1$, i.e., 
$A= e\, B + \bar e\, (B^t)^{-1}$ with $B\in {\rm SL}_n(\bbr)$.  
Thus ${\rm SU}_n(\bbd) \cong {\rm SL}_n(\bbr)$.  Note that $(B^t)^{-1} = (\det B)^{-1}\wt B$.

 \vskip.3in

\centerline{\headfont \  \D.\  The Special Lagrangian Calibration. }  
\medskip

The Special Lagrangian calibration, and its associated differential equation, which were introduced in [HL$_1$],
have a complete analogue in $\dn$.  This work, due to Jack Mealy [M$_{1,2}$], will be presented in this
and the following sections.

\medskip
\noindent
{\bf Some Preliminaries.}
A real $n$-plane $P$ in $\dn$ is {\bf Lagrangian} if $\o\bigr|_P=0$.  
Equivalently, $\bra {z'}{Tz} = \o(z',z)=0$ for all $z,z'\in P$, i.e., $z\in P \ \Rightarrow
\ Tz \in P^\perp$.  Let $\LAG$ denote the set of all Lagrangian $n$-planes.  A real $n$-plane
$P$ in $\dn$ is {\bf space-like} if the inner product $\bra \cdot\cdot $ on $\dn$, when restricted to 
$P$, is positive definite. In particular, $P\cap \tau\rn = \{0\}$ so that $P$ can be graphed over
$\rn$ in $\dn=\rn\oplus\tau\rn$ by a matrix $A\in M_n(\bbr)$. Pulling back the inner product
 $\bra \cdot\cdot $ on $\dn$ to $\rn$ by the graphing parameterization
  $x\mapsto x+\tau Ax$ yields $\bra xx-\bra {Ax}{Ax}$. Hence 
  $$
  P\ \ {\rm is \  spacelike}\qquad\iff\qquad I-A^tA\ >\ 0.
 \eqno{(\D.0)}
$$ 

Let $\Gs$ denote the Grassmannian of oriented space-like $n$-planes in $\dn$.
There are two connected components $\Gsp$ and $\Gsn$. The graph of $A$, when
given the orientation from $\rn$, is defined to be in $\Gsp$, and when given the opposite 
orientation, is defined to be in $\Gsn$.

Recall that an oriented $p$-plane in a real vector space $V$ does {\sl not} uniquely determine
a simple (or decomposable) vector $\x\in\L^p V$ by choosing an oriented basis $\e_1,...,\e_p$ and 
defining $\x=\e_1\wee\e_p$.  However, $\x$ is determined up to a positive scale.
If $V$ is equipped with a positive definite inner product, then  $\x\in\L^p V$  is uniquely determined
by requiring $\e_1,...,\e_p$ to be an oriented orthonormal basis, and the Grassmannian of oriented
$p$-planes  in $V$, denoted $G(p,V)$, can be identified with a subset of $\L^p V$.
This is not possible if $V$ is equipped with an inner product such as the one on $\dn$.  
However, if $\xi\in\Gs$ is an oriented space-like $n$-plane in $\dn$, then again by requiring the 
basis $\e_1,...,\e_n$ to be orthonormal with respect to the induced inner product, this
subset $\Gs$ of the full Grassmannian can (and will be) identified with a subset of $\L^n_\bbr\dn$.

Finally, denote by $\LAGt$ the Lagrangian $n$-planes in $\dn$ which are both oriented and space-like.
The component $\LAGt^+$ of $\LAGt$ containing the standard oriented
$\overrightarrow{\bbr}^n$ in $\dn=\rn\oplus \tau\rn$ will be referred to as the set 
of {\sl positive oriented space-like Lagrangian $n$-planes}.  The sets $\LAGt^+\ss\LAGt\ss \Gs$ 
are all subsets of $\L^n_\bbr\dn$.

\medskip
\noindent
{\bf The Special Lagrangian Inequality.}  As before, let $dz= dz_1\wee dz_n$ denote the {\sl standard}
$(n,0)$-form on $\dn$.  

\Lemma{\D.1} {\sl
Suppose $\x' \in \L^n_\bbr\dn \cong \L^n \bbr^{2n}$ is a real $n$-vector and $A\in M_n(\bbd)$.
Then
}
$$
(dz)(A\x') \ =\ \left( \det_\bbd A    \right) (dz)(\x').
$$
\pf
Note that $(dz)(A\x')  =  \left( A^t (dz)\right) (\x') = \left(\det_\bbd A^t    \right) (dz)(\x')
= \left(\det_\bbd A    \right) (dz)(\x')$.\qed\medskip

In analogy with  the complex case each $\x\in\LAGt^+$ has a {\sl phase}
$\theta\in\bbr$.

\Cor{\D.2} {\sl
If $\x\in\LAGt^+$, then}
$$
(dz)(\x) \ =\ e^{\tau\theta}\ \ \ {\sl for \ some\ } \theta\in \bbr
$$
\pf
Since $(dz)(\x_0) =1$ where $\x_0 = { \oa \bbr}^n$, if $A\in {\rm U}^+_n(\bbd)$, then
$$
(dz)(A\x_0) \ =\ \det_\bbd A
\eqno{(\D.1)}
$$
By (3.4)(3), each $\x\in\LAGt^+$ is of the form $\x=A\x_0$ with 
$A\in {\rm U}_n^+(\bbd)$.\qed
\medskip

Since $A\in {\rm U}_n^+(\bbd)$ fixes $\overrightarrow{\bbr}^n
 \ \iff\ A\in {\rm SO}_n(\bbr)\ss M_n(\bbr)$, this proves that
$$
\LAGt^+ \ =\ {\rm U}_n^+(\bbd)/{\rm SO}_n(\bbr).
\eqno{(\D.2)}
$$ 

\Def{\D.3}  A positively  oriented space-like Lagrangian $n$-plane $\x \in \LAGt^+$ is said
to be {\bf (split) special Lagrangian} or to have {\bf phase zero} if $\x=A\x_0$ for some 
$A\in  {\rm SU}_n(\bbd)$.  The space of all (split) special Lagrangian $n$-planes in $\dn$ is denoted by
$\SLAGt$.

Note that 
$$
\SLAGt \ =\ {\rm SU}_n(\bbd)/{\rm SO}_n(\bbr).
\eqno{(\D.3)}
$$ 
since ${\rm SO}_n(\bbr)$ is contained in ${\rm SU}_n(\bbd)$,
and   (\D.2)  holds.

\Theorem{\D.4. (Mealy)}
{\sl
$$
({\rm Re}\, dz)(\x) \ \geq\  1 \fa \x\in \Gsp
$$
with equality if and only if $\x\in\SLAGt$.}

\pf
Let $\be=(1,0)$ and $T\be=(0,1)$ denote the standard basis for $\bbd\cong \bbr^2$.
Then the dual basis for $(\bbr^2)^*$ is $\be^* = dx$ and $(T\be)^* = dy$.
Employing this same notation in $\dn$, we have
$$
dz\ =\ \left( \be_1^* + \tau (T\be_1)^*   \right)
\wee
 \left( \be_n^* + \tau (T\be_n)^*   \right).
\eqno{(\D.4)}
$$
Expanding out this product expresses $dz$ as a $2^n$-fold sum.
Each term $\a$ in the sum is the product of $n$ of the $2n$ axis covectors
$\be_1^*,  \tau (T\be_1)^*,...,\be_n^*, \tau (T\be_n)^*$.  However, for any $j=1,...,n$,
$$
{\rm both}\ \ \be_j^*  \ \ {\rm and\ \ }  \tau (T\be_j)^*\ \ {\rm cannot\ be\ factors\  of\ }\ \a.
\eqno{(\D.5)}
$$

By Proposition A.3 each $\x\in \Gsp$ is unitarily equivalent to 
$$
\x'\ =\      \be_1  \w      \left( \cosh \theta_1 \be_2^* + \sinh \theta_1  (T\be_1)^*   \right) \w
 \be_3  \w      \left( \cosh \theta_2 \be_4^* + \sinh \theta_2  (T\be_3)^*   \right) \w
\cdots,
\eqno{(\D.6)}
$$
that is, $\x=A\x'$ with $A\in {\rm U}^+_n(\bbd)$.
By (\D.4), (\D.5) and (\D.6) we have
 $$
 (dz)(\x') \ =\  \cosh \theta_1 \cdots \cosh \theta_{[{n\over 2}]}.
 $$
Hence by Lemma \D.1
$$
 (dz)(\x) \ =\   e^{\tau\theta} \cosh \theta_1 \cdots \cosh \theta_{[{n\over 2}]},
\eqno{(\D.7)}
$$
where $e^{\tau\theta} = \det_\bbd A$.  In particular
$$
 ({\rm Re}\,dz)(\x) \ =\   \cosh \theta \cosh \theta_1 \cdots \cosh \theta_{[{n\over 2}]}.
\eqno{(\D.8)}
$$
Hence $({\rm Re}) (\x) \geq 1$ and $=1$ if and only if all the angles are zero, i.e., 
 $\x' = \be_1\wee \be_n = \x_0$, $\x = A\x' = A\x_0$, and $\det_\bbd A=1$.
\qed

\Cor{4.5} {\sl
An oriented real $n$-plane $\x$ in $\dn$ is (split) special Lagrangian if and only if 
$\x\in \LAGt^+$ and   ${\rm Im} \,dz\bigr|_\x =0$.
}

\bigskip
\noindent
{\bf The Null Viewpoint.}
We now revisit the material above in null coordinates.
Several points are worth mentioning.  Note first that (4.3) naturally becomes
$$
\SLAGt \ =\ {\rm SL}_n(\bbr)/{\rm SO}_n(\bbr).
\eqno{(\D.3)'}
$$ 
since $A\in {\rm SU}_n(\bbd)$ if and only if $A=e\,B + \bar e\, (B^t)^{-1} = e \, B+\bar e\, \wt B$
for some $B\in {\rm SL}_n(\bbr)$. 

Note also that
$$
\qquad\qquad
dz\ =\ \half(du+dv) + \smfrac \tau 2 (dv - du), \qquad{\rm or \  }
\eqno{(\D.9)}
$$ 
$$
{\rm Re}\,dz\ =\ \half(du+dv)  \and  
{\rm Im}\,dz\ =\  \half (dv - du).
\eqno{(\D.10)}
$$ 
To prove this note that $dz_j = e\,du_j+\bar e\,dv_j$ implies that
$dz =  e\,du+\bar e\,dv$.

Now Corollary \D.5 has the null form:

\Cor{4.5$'$} {\sl
An oriented real $n$-plane $\x$ in $\dn$ is (split) special Lagrangian if and only if 
$\x\in \LAGt^+$ and  $du\bigr|_\x = dv\bigr|_\x$.
}

   \vskip .3in

\centerline{\headfont \ \SS.\  Split SLAG Submanifolds. }  
\medskip

The algebraic calculations of the previous sections apply to submanifolds.
In what follows all submanifolds are assumed to have dimension $n$.

\Def{\SS.1}  A closed oriented   $C^1$-submanifold $M$ of an 
an open subset of $\dn$ is a {\bf split SLAG} (or {$\wt{\bf SLAG}$)  {\bf submanifold} 
 if the oriented tangent space $\oa {T_z}M\in \SLAGt$ 
for all $z\in M$.

\Theorem{\SS.2} {\sl
A closed oriented   $C^1$-submanifold $M$ of an 
an open subset of $\dn$ is {split SLAG} if and only if 
\medskip

(1)\ \ $M$ is Lagrangian, i.e., $\o\bigr|_M=0$,
\medskip

(2)\ \ $\psi\bigr|_M=0$ where $\psi \equiv {\rm Im} \,dz = \half(dv-du)$, and
\medskip

(3)\ \ $M$ is positive space-like.
}

\pf Apply Corollary \D.5 and (\D.10).\qed
\medskip

If $M$ satisfies only conditions (1) and (2),  it will be called an {\sl unconstrained}
 (or {\sl not necessarily space-like}) split SLAG submanifold.

Each oriented space-like submanifold $N$ inherits a volume form vol$_N$ from $\dn$.

\Theorem{\SS.3} {\sl
Suppose that $(M,\partial M)$ is a compact  oriented submanifold with boundary in $\dn$.
If $M$ is split SLAG, then $M$ is volume-maximizing, i.e., 
$$
{\rm vol}(M)\ \geq\ {\rm vol}(N)
\eqno{(\SS.1)}
$$
for any other positive space-like compact submanifold $N$ with $\partial N=\partial M$.

Equality holds in (\SS.1) if and only if $N$ is also split SLAG.
}
\pf
We have 
$$
{\rm vol}(M)\  =\  \int_M {\rm Re} \, dz  \  =\  \int_N {\rm Re} \, dz \  \geq \  {\rm vol}(N)
$$
by Theorem  \D.4, Stokes' Theorem, and Theorem \D.4 again.
Equality in (\SS.1) implies $ {\rm Re} \, dz\bigr|_N  \equiv  {\rm vol}(N)$ 
and therefore $N$ is split SLAG by  Theorem \D.4 once again.
\qed

\Remark {\SS.4. (Split SLAG$_\theta$-Calibrations)}
As in the complex case, each rotation of the 
calibration ${\rm Re} \, dz$ above gives a new calibration ${\rm Re}  (e^{-\tau \theta} dz)$
on positively oriented space-like submanifolds. The corresponding calibrated submanifolds
are the Lagrangian submanifolds with the property that the restriction 
$$
{\rm Im} (e^{-\tau \theta} dz)\bigr|_M\ =\ 0
$$
or equivalently, that ${\rm Re}  (e^{-\tau \theta} dz))\bigr|_M = d{\rm vol}_M$.

Note that  Corollary \D.2 provides a fibration
$$
\LAGt^+\ \arr \ \bbr
$$
whose  fibre at $\theta\in \bbr$ 
is exactly the set 
$
\SLAGt_{\theta}
$
of   $n$-planes in $\Gsp$  calibrated by ${\rm Re}  (e^{-\tau \theta} dz)$.

We shall not discuss these other cases since they are completely equivalent to the 
case $\theta=0$.  If $A \in {\rm U}^+_n(\bbd)$ with $\det_\bbd(A)= e^{\tau\theta}$,
then  this isometry $A:\dn\to\dn$ carries $\SLAGt$-submanifolds to 
$\SLAGt_{\theta}$-submanifolds.

   \vskip .3in

\centerline{\headfont \ \DD.\  Split SLAG Graphs in $\dn$. }  
\medskip

Locally the graph of a smooth map $F:\rn\to\rn$ is Lagrangian in $\bbr^{2n}$ if and only if
$F=df$ for some potential function $f$.  This is because $\o=\sum_j dx_j\w dy_j$  vanishes when
restricted to $M$ if and only if  $\sum_j  y_j dx_j \bigr|_M$ is $d$-closed.  This is equivalent to 
requiring the Jacobian of $F$ to be symmetric, in which case it equals the hessian of the potential 
$f$.

Any space-like submanifold $M$ is locally the graph of a smooth map from $\rn$ to $\tau\rn$
since $T_zM \cap \tau\rn = \{0\}$ for all $z\in M$. Furthermore, in 
Proposition \R.6 below we show that if $M$ is a  closed space-like Lagrangian submanifold
of a tube $\O\times \tau \rn \ss \rn\times \tau\rn = \dn$, then the projection $\pi\bigr|_M:M \to\O$
is a covering map.  In particular, if $\O$ is simply-connected, then each connected component
of $M$ is the graph of a gradient $df$ for some potential function $f:\O\to \bbr$.
This makes the following theorem particularly relevant.

\Theorem {\DD.1. (Mealy)}  {\sl
Suppose that $f$ is a smooth real-valued function defined on an open subset
of $\rn$.  Let $M$ denote the oriented graph of the gradient 
$F=df$ in $\dn = \rn \oplus \tau \rn$.  Then $M$ is split special Lagrangian 
\ $\iff$\  $M$ is space-like and
$$
\sum_{k=0}^{[(n-1)/2]} \s_{2k+1}(\Hess\, f)\ =\ {\rm Im} \, \det_\bbd\left( I+\tau \Hess\, f\right) \ =\ 0.
\eqno{(\DD.1)}
$$
Furthermore $M$ is space-like\ $\iff$\ 
$$
-I\ <\ \Hess\,f \ <\ I.
\eqno{(\DD.2)}
$$
or, equivalently, the graphing map of $df$ is uniformly Lipschitz with Lipschitz constant 1.}

\pf
The tangent plane $P$ to $M$ at a point is the graph of $A\in   M_n(\bbr)$ in $\rn\oplus\tau\rn$,
and is parameterized by $x\mapsto x+\tau Ax$,  $x\in\rn$.  Pulling back the symplectic form $\o$ from 
$P$ to $\rn$ gives $\o((x,Ax), (x',Ax')) = \bra x{Ax'} - \bra{Ax}{x'}$ so that, as noted above, 
$P$ is Lagrangian \ $\iff$\ $A$ is symmetric.  Pulling back the quadratic form $\bra \cdot\cdot$
yields $\bra x x - \bra {Ax}{Ax}$.  Hence, as noted in (\D.0),
$$
 P \ \equiv\ {\rm graph}\, A \ {\rm is\ spacelike} \ \ \iff\ \ A^t A\ <\ I.
\eqno{(\DD.2)'}
$$

It suffices to prove Theorem \DD.1 when $M$ is a plane.  This case is handled explicitly
in the following.

\Prop{\DD.2}{\sl
Suppose $\x={\rm graph}\,A$ is the oriented graph in $\dn=\rn\oplus \tau\rn$ of a matrix $A\in M_n(\bbr)$.
Then $\x\in \SLAGt$ \ $\iff$\ $A$ is symmetric with $-I<A<I$ and  
}
$$
{\rm Im} \, \det_\bbd\left( I+\tau A \right) \ =\ 0.
\eqno{(\DD.3)}
$$
\pf
We have already shown that $\x$ is space-like and Lagrangian (i.e., $A\in\LAGt$)  \ $\iff$\ $A$ is symmetric and
$A^2<I$, (which is equivalent to $-I< A <I$ for symmetric $A$).  Since the space of such oriented
graphs is connected and contains $\x_0 = {\oa\bbr}^n$ (take $A=0$), this proves that
$$
 \x = {\rm graph}\,A\ \in \ \LAGt^+ \quad\iff\quad A\ \ {\rm is\ symmetric\ and} \ -I < A <I.
\eqno{(\DD.4)}
$$

Now the pull-back of $dz$ under the $\bbd$-linear map $I+\tau A :\dn\to \dn$  is equal to
$\det_\bbd(I+\tau A) \,dz$.  Hence, the pull-back of $dz\bigr|_\x$ to $\rn$ under the
graphing map $x\mapsto x+\tau Ax$ equals $\det_\bbd(I+\tau A) \, dx$.
Hence (\DD.3) holds \ $\iff$\  Im$\, dz\bigr|_\xi =0$.

Corollary \D.5 completes the proof of Proposition \DD.2 and Theorem \DD.1.
\qedqed

\bigskip
\noindent
{\bf The Null Viewpoint.}  
It is frequently more interesting to consider graphs over the null $n$-plane $e\,\rn$ in
$\dn= e\, \rn \oplus\bar e\, \rn$.

Each space-like submanifold $M$ in $\bbd^n$ is locally the graph of a smooth map 
from  a domain $\O\sse\rn$ to $\bar e\rn$ since $T_zM \cap \bar e\rm = \{0\}$ for all $z\in M$.
Furthermore, since  $\o=\sum_j du_j\w dv_j$  we see that  if  $M$ is also Lagrangian 
and $\O$ is simply-connected, then 
this graph is the gradient $df$ of a potential function $f$ on $\O$.

\Theorem{\DD.3}  {\sl
Suppose that $f$ is a smooth function defined on an open subset of $\rn$.
Let $M$ denote the oriented graph of the gradient $F=df$ in $e\rn\oplus\bar e \rn \equiv \dn$.
Then
\bigskip
\centerline{ $M$ is split special Lagrangian \ $\iff$\  $f$ is convex and $\det_\bbr \Hess\, f\equiv 1$.}
}
\medskip
\pf
 The standard plane
$\x_0 = \overrightarrow{\bbr}^n$  is the graph of $A=I \in M_n(\bbr)$.
For any $A\in M_n(\bbr)$,  the previous discussion
shows that $\x={\rm graph}\, A$ is Lagrangian \ \ $\iff$\ \  $A$ is symmetric.
The quadratic form $\bra \cdot\cdot$ restricted to the graph of $A$ pulls back to
$\bra {e u+\bar e Au}{e u+\bar e Au} = 2\bra{e u}{\bar e Au}\ =\ \bra u{Au}$ since
$eu$ and $\bar e Au$ are null and $\bra e{\bar e} = \half$ in $\bbd$.  Thus:
$$
{\rm graph} \, A\ \ {\rm is\ spacelike} \quad\iff\quad A+A^t  \ >\ 0.
\eqno{(\DD.5)}
$$
 By a connectivity argument it follows that: for   the oriented graph $\x$ of $A\in\M_n(\bbr)$
 $$
 \x = {\rm graph}\,A \in \LAGt^+\quad\iff\quad A\ {\rm is\ symmetric\ and\ } A\ >\ 0.
\eqno{(\DD.6)}
$$

\Prop{\DD.4}  {\sl
One has $\x\in\SLAGt\ \ \iff\ \ \x$ is the oriented graph 
over $e\,\rn$ in $\dn= e\,\rn \oplus\bar e\, \rn$ of a positive definite symmetric 
matrix $A\in M_n(\bbr)$ satisfying
}
$$
\det_\bbr \, A\ =\ 1.
\eqno{(\DD.7)}
$$
\pf
As noted above,  each space-like $n$-plane $P$ can be graphed over $e\,\rn$ since $P\cap e\rn = \{0\}$.
Under the graphing map $u\mapsto eu+\bar e Au$, the $n$-form $du$ pulls back to $du$ and 
the $n$-form $dv$ pulls back to $(\det_\bbr A)du$.  Thus by Corollary \D.5$'$ we have
$\x = {\rm graph} \,A \in \SLAGt\  \iff \ \x\in \LAGt^+$ and $\det \, A=1$.
This completes the proof of both Proposition \DD.4 and Theorem \DD.3.\qedqed

\vskip.3in


\centerline{\headfont \ \R.\  Removable Singularities for   split  SLAG  varieties. }  
\medskip

In this section we show that  split SLAG subvarieties tend to be less singular than their
elliptic cousins.  For example, any connected  split  SLAG submanifold in a 
simply-connected tube domain
$\O\times\tau \rn$ is the graph of a mapping with potential function satisfying (\DD.1).  
Moreover, singularities in
codimension $> 2$ tend to be removable.  In particular, there are no irreducible
 split  SLAG cones with isolated singularity in $\bbd^n$ for $n\geq 3$.

The proofs of these results, which follow, illustrate an interesting interplay between the
$z$-coordinate and the null-coordinate viewpoints.  This shift of coordinates, sometimes
called the {\sl Cayley transformation}, turns out to be quite useful in the analysis.

The first main theorem of this section is the following.

\Theorem{\R.1. (Removable Singularities I)} {\sl
Let $\O\ss\rn$ be a convex domain in $\rn$ and consider the 
``tube'' domain $\Omega\times \tau\rn \ss  \rn\oplus \tau \rn=\dn$.
Let $\Sigma\ss\O$ be a compact subset of Hausdorff $(n-2)$-measure zero. Then any closed, 
connected  split  SLAG submanifold 
$$
M\ss (\O-\Sigma)\times \tau\rn
$$ 
has closure in $\O\times \tau\rn$ which is the graph of a real analytic mapping
$F:\O\to\rn$.  Furthermore, $F=\nabla f$ where $f:\O\to \bbr$ satisfies the equation
(\DD.1).
}

 \Cor{\R.2. (Absence of Cones)} {\sl 
 Suppose $C\ss \bbd^n$  is a    split  SLAG  cone which is regular and connected outside the origin.
 Then $C$ is a $\SLAGt$ $n$-plane.
 }

\medskip
\noindent{\bf Proof of Corollary \R.2.}  The proof is immediate for $n\geq 3$ -- one takes $\O=\rn$ and 
$\Sigma = \{0\}$. The case $n=2$ follows from Theorem F.3 in Appendix F and the fact that 
complex cones of complex dimension 1 are planes.\qed

\Remark{\R.3}  This is in sharp contrast to the non-split case where ${\rm SLAG}$-cones exist
in all dimensions and there are many ${\rm SLAG}$-varieties with singularities of high codimension.
(See [HL$_1$],  [J$_*$],[H$_*$], [HK$_*$]  and the references therein). Also in coassociative geometry there
is a coassociative cone, smooth outside 0 in $\bbr^7$ which is the graph of a Lipschitz map [HL$_1$].

\Remark{\R.4} The codimension-2 hypothesis in Theorem \R.1 is close to best possible.
There do exist   split  SLAG  subvarieties in all dimensions $n\geq2$ with codimension-2
singularities.  When $n=2$ they appear because     split  SLAG  subvarieties are
 complex curves for a certain complex structure (see Appendix E). The higher dimensional
cases are just products of these with $\bbr^{n-2}$.

On the other hand we still do not know of any  non-trivial (irreducible)   split  SLAG  cones.

 \Cor{\R.5. (Removable Singularities II)} {\sl 
Let $\O\ss\rn$ be a  domain in $\rn$ and suppose
 $\Sigma\ss\O$ is any closed subset of linear-measure zero.
 Then any closed, 
   split  SLAG  submanifold 
$
M\ss (\O-\Sigma)\times \tau\rn
$
has closure in $\O\times \tau\rn$ which is an immersed 
real analytic   split  SLAG  submanifold of $\O\times \rn$.
}
\medskip\noindent
{\bf Proof of Corollary \R.5.} Fix $(x,y)\in \overline M$. Choose a ball $B=B_r(x)\ss\O$ of radius
$r$ such that $\partial B \cap \Sigma = \emptyset$.  One now applies Theorem \R.1
to the tube  domain $B \times \rn$ to conclude that the closure of each connected component
of $M\cap (B-\Sigma) \times \rn$  is the graph of an analytic mapping $B\to\rn$.  \qed

\medskip
\noindent{\bf Proof of Theorem \R.1.}  We begin with the following observation.
 
\Prop{\R.6} {\sl
Let $\O\ss\rn$ be a domain in $\rn$ and consider the 
``tube'' domain $T\equiv \Omega\times \tau\rn \ss \rn\oplus \rn=\dn$ with projection $\pi:T\to \O$.
 If  $M \ss T$ is a closed (embedded) submanifold which is space-like and Lagrangian, 
 then $\pi\bigr|_M : M\to\O$ is a covering map.

 If, in addition, $\O$ is simply-connected, then each component of $M$
  is the graph of $\nabla f$ for some
 potential function $f:\O\to\bbr$, and this graph mapping $F=\nabla f$ is 1-Lipschitz.
 }
\pf
To begin we note that
from (\DD.2) we have the following:
$$
\| \pi_*(v)\| \ \geq\ {1\over \sqrt 2} \|v\| \fa v\in T(M).
\eqno{(\R.1)}
$$
Now fix a point $x\in \O$ such that $x=\pi(y)$ for some $y\in M$.  Let
$\g:[0,1] \to \O$ be a smooth  curve with $\g(0)=x$.  From \R.1 one easily concludes
that there exists a lift $\wt \g :[0,1]\to M$ with $\wt\g(0)=y$ and
$\pi\circ \wt\g= \g$.  (The set of $t\in [0,1]$ such that $\g\bigr|_{[0,t]}$ can be lifted, is both open and closed.)

Fix a closed ball $B= B(x,r)$ about $x$ which is contained in $\O$.
Then lifting radial lines from $x$ (as above) gives a mapping $\vf : B\to M$ with 
$\vf(x)=y$ and $\pi\circ \vf= {\rm Id}_B$.  It now follows that for each connected 
component $M_0$ of $M$, the mapping 
$\pi_0=\pi\bigr|_{M_0}:M_0\to \O$ is a covering map, i.e., $\pi_0$
is  surjective and every point $x$ has a neighborhood $B$ which is evenly covered
by $\pi_0$. The first assertion is proved.

 If $\O$ is simply-connected and $M$ connected, then by elementary covering space theory 
$\pi=\pi\bigr|_{M}:M\to \O$ is a diffeomorphism.  Hence $M$ is the graph of a function
$F:\O\to \rn$, which, since $M$ is Lagrangian and $\O$ is simply-connected,
is the gradient $F=\nabla f$ of a scalar function $f:\O\to\bbr$.

Since $M$ is also space-like, we have from (\DD.2) that
$$
-I \ <\ {\rm Jac} (F)\ <\ I \qquad {\rm on\ } \O
\eqno{(\R.2)}
$$
Hence $F$ is 1-Lipschitz and the proof is complete.\qed

\medskip

Returning to  Theorem \R.1, we apply Proposition \R.6
 to the simply-connected domain $\O-\Sigma$
to conclude that $M$ is the graph of a gradient  $F=\nabla f:\O-\Sigma\to \rn$
which satisfies (\R.2). The map $F$ extends uniquely to a 1-Lipschitz mapping
$\hat F:\O\to \rn$ which can be written as the gradient $\hat F=\nabla \hat f$ of a unique extension 
of $f$ to $\O$.  It follows that the closure of $M$ in $\O$ is
$$
\overline M  \ =\ M\cup \Sigma_M   
$$
where
$$
\Sigma_M \equiv (I\times \hat F)(\Sigma) \ \ {\rm is\ compact \  with \ Hausdorff }\ (n-2) \ {\rm measure\ zero.}
\eqno{(\R.3)}
$$

We now pass to null coordinates $(u,v)$ on $\dn$ and consider the projection of $\overline M$
onto the $u$-axis. Composing this projection with the graphing map $F$ gives a mapping
$G=F+I:\O\to \rn$ (i.e.,  $G(x)=F(x)+x$ for $x\in \O$).  
Note that $G=\nabla g$ where $g(x)= f(x) + \half |x|^2$.
From  (\R.2) we see that this map satisfies
$$
0\ <\ {\rm Jac\, G} \ <\ 2I \qquad {\rm on\ } \O-\Sigma.
$$
In particular, the potential 
$g$ is strictly convex on $\O-\Sigma$.
It follows that for any two points $x, y\in \O$, we have that
$$
G(x)\ \neq\ G(y) \ \ {\rm  unless\  the\  segment \ }  \overline {xy} \ss \Sigma.
\eqno{(\R.4)}
$$
In particular, $G$ is injective on $\O-\Sigma$.  Thus $G\bigr|_{\O-\Sigma}$ is a 
diffeomorphism onto its image. Moreover we have the following.

\Lemma{\R.7}  {\sl The map $G\bigr|_{\O-\Sigma} : \O-\Sigma\to \O'-\Sigma'$
is a diffeomorphism, where $\O'\equiv G(\O)$ and $\Sigma'\equiv G(\Sigma) \ss\ss\O'$.
}

\pf
By (\R.4) and the strict convexity we know that if $x\in \O-\Sigma$ and $y\in \Sigma$, then
$G(x)\neq G(y)$.  Hence, $G(\O-\Sigma)\cap G(\Sigma)=\emptyset$.
By shaving $\O$ slightly we may assume it has a smooth boundary $\partial \O$ which
does not meet $\Sigma$ and along which $G$ is a diffeomorphism. Evidently $G(\Sigma)$
will be a compact subset of the bounded component of $\rn-G(\bo)$. \qed

\medskip

Since $G$ is 2-Lipschitz, the Hausdorff 2-measure $\ch(\Sigma')=0$. It follows that
$\O'-\Sigma'$ is simply-connected.  Now note that $\overline M$ is the graph in null coordinates 
of a mapping $v=H(u)$ where  
$H:\Omega' \to \rn$ is well defined and smooth except possibly over  points of $\Sigma'$. (Points 
where it is not well defined correspond to 
points lying on non-trivial straight line segments in $\Sigma$). Hence, in $\O'-\Sigma'$
we know that $H=\nabla h$
where $h$ is a smooth strictly convex function which satisfies the equation
$$
\det\, \Hess(h) \ =\ 1
\eqno{(\R.5)}
$$
on $\O'-\Sigma'$.  Here we are using the assumption that $M$ is $\SLAGt$ and Theorem \DD.3.

Now any  convex function defined on $\O'-S$ where $S\ss\O$ is closed with $\ch^{n-2}(S)=0$,
 extends as a convex function across $S$.
Thus $h$ is well defined and convex 
on $\O'$ 

We claim that $h$ is a weak solution of equation (\R.5)  (and therefore a viscosity solution
by [C$_3$, Lemma 3])  on all of $\O'$. To see this note that for any Borel
set  $\co\ss\O'-\Sigma'$ equation (\R.5) implies that
$$
{\rm meas} (\nabla h (\co) )\ =\ {\rm meas}  (\co).
$$
However, by (\R.3)   this holds also for any Borel set $\co\ss\O'$.

We now apply  deep results of Caffarelli (Theorem  \R.8 below)  to conclude 
that $h$ is real analytic in $\O'$.  Hence, $M = {\rm graph}(\nabla h)$ is a real analytic 
submanifold of $\O\times \rn$. Thus the graph of $F=\nabla f$ over $\O$ is real analytic.
Since $F$ is 1-Lipschitz, it is also real analytic on $\O$.\qed
\medskip

The following is well-known but not explicitly stated in the literature.
\Theorem{\R.8.  (Caffarelli [C$_{2,3}$  ])}  {\sl
Let $h$ be a viscosity solution of (\R.5) on a bounded domain $\O\ss\rn$ with smooth
boundary $\partial \O$.  Assume $h$ is smooth in a neighborhood of $\partial \O$
in $\overline \O$.  Then $h$ is real analytic in $\O$.}

\medskip
\def\bdu{{\bf u}}\def\bdv{{\bf v}}

\Ex{\R.9} Take the union of the standard $\rn\ss\dn$ with the $n$-plane: 
$y_1= \cdots = y_{n-2} = 0, \ y_{n-1}=\l x_{n-1}, \ y_n= - \l x_n$.  Both are split SLAG. This union is singular 
(but not irreducible) along the intersection $(\bbr^{n-2} \times \{0\})\times \{0\}$.

\Ex{\R.10}  The following classical example of Pogorelov is of interest in this context.
Let  $(u,v)$ be null coordinates on $\dn$, 
and write $u=(u_1,\bdu)$ and $v=(v_1,\bdv)$  for  $\rn=\bbr\times \bbr^{n-1}$.
Consider the potential
$$
\vf(u_1, \bdv) \ =\ -k|\bdv|^{2(n-1)\over n-2} f(u_1)
$$
where $f(t)$ satisfies $f''+f=0$ and $f>0$.  The graph $M$ of the gradient of $\vf$ is:
$$
\eqalign
{
v_1\ &=\ -k|\bdv|^{2(n-1)\over n-2} f'(u_1)  \cr
\bdu \ &=\ -k  \smfrac{2(n-1)}{n-2}  |\bdv|^{2\over n-2} f(u_1) \bdv.   
}
$$
Note that $M$ is real analytic if $n=3$, but just $C^{1,\a}$ along the $u_1$-axis if $n\geq 4$
(where $\a = ([{n\over 2}]-\half)^{-1}$ when $n$ is odd, and 
$\a = ({n\over 2}- 1)^{-1}$ when $n$ is even).
The tangent plane at points on the $u_1$-axis is the $(u_1, \bdv)$-plane,
so $M$ is being graphed over its tangent plane at these points.

With the right choice of $k$, $M$ is an example of a split SLAG manifold outside of a 
singular line which is null.  In particular, when $n=3$, $M$ is a real analytic unconstrained split SLAG
submanifold.  Outside the $u_1=x_1-y_1$ axis it is  split SLAG, but  at points on this axis
it is not split SLAG since this 
axis is tangent to $M$ and null.

\Remark{\R.10} Although the discussion in this section was carried out without defining the notion
of a (singular) split SLAG variety, we call Appendix D to the reader's attention for this definition.

  \vskip.3in

\centerline{\headfont \  \E. \  Double Manifolds (D-Manifolds). }  
\medskip

Certain standard notions from complex analysis can now be carried over to $\bbd$-analysis.

\Def{ \E.1}  A {\sl double manifold } or {\sl $\bbd$-manifold} is a smooth manifold $X$ equipped with an
atlas of charts $\{(U_\a,\psi_\a)\}_\a$, $\psi_\a:U_\a\to \bbd^n$, whose transition  functions
are $\bbd$-holomorphic and orientation-preserving, that is, the transition functions are
of the form
$$
z'\ =\ F(z)\qquad {\rm with\ \ } {\partial F\over \partial \bar z}\ =\ 0\and
 \det_\bbd\left( {\partial F\over \partial  z}\right)\ \in \ \bbd^+
\eqno{( \E.1)}
$$
or equivalently
$$
u'\ =\ f(u), \ \ v'\ =\ g(v) \qquad {\rm with\ \ } f, g \ \ {\rm orientation \ preserving}.
\eqno{( \E.2)}
$$
  
Let's examine these manifolds from the two points of view.  
In the complex picture, we have a Dolbeault decomposition 
$$
\left( \L^k T^*X \right ) \otimes_\bbr \bbd \ =\ \bigoplus_{p+q=k} \L^{p,q}(X)
$$
where, in local holomorphic coordinates $(z_1,...,z_n)$,   $\L^{p,q}$ is spanned by 
forms of the type $dz_{i_1}\wedge\cdots \wedge dz_{i_p}\wedge d\bar z_{j_1}\wedge\cdots \wedge d\bar z_{j_q}$.
Let $\ce^{p,q}(X)$ denote the smooth sectons of $\L^{p,q}(X)$.  Then the   operators
$$
\partial : \ce^{p,q}(X)\ \arr\ \ce^{p+1,q}(X)
\and
\dbar : \ce^{p,q}(X)\ \arr\ \ce^{p,q+1}(X),
$$
already given in local holomorphic coordinates by (3.1),
are well defined on the manifold $X$.
A $\bbd$-valued function $F$ defined on an open subset of $X$ is holomorphic if and only if $\dbar F=0$.

Suppose now that we rewrite our local charts on $X$ in terms of null coordinates 
$(u,v)$ in $\dn = e\cdot \bbr^n\oplus \bar e\cdot \bbr^n$.   Then from (\E.2) it is clear that $X$ is furnished with
two transversal $n$-dimensional  oriented foliations $\cf^+$ and $\cf^-$ given locally by the
$u$ and $v$ coordinate planes respectively. Taking the tangent plane fields $N^\pm  \equiv \cf^\pm$ to be the 
$\pm1$-eigenspaces defines the {\bf structure endomorphism}
$$
\bbt : T(X) \ \to\ T(X)  \qquad {\rm with\ \ } \bbt^2\ =\ Id.
\eqno{(\E.3)}
$$
which represents scalar multiplication by $\tau$ on the tangent spaces (cf. (2.3)).
(This is the analogue of $J : T(X) \to T(X)$ with $J^2=-Id$  in the complex case.)
Indeed from this point of view a $\bbd$-manifold is simply a pair $(X,\bbt)$, where $\bbt^2=Id$
with eigenbundles $N^\pm$ of equal dimension both of which are integrable. (If one 
drops the integrability condition, this becomes an {\sl almost $\bbd$-manifold}.)
Even more simply, a $\bbd$-manifold is a 2$n$-dimensional manifold with a pair of 
transversal $n$-dimensional oriented foliations.
When $n=1$ this is equivalent to a (fully oriented) conformal Lorentzian structure on the surface.

\Remark{ \E.2} Any product $X=M_1\times M_2$ of $n$ manifolds is certainly a $\bbd$-manifold 
(with the foliations given by the factors). 
However, $\bbd$-manifolds can be much more complicated. 
Let $\cs_k$ be a 2-dimension foliation on a compact 3-manifold $M_k$, $k=1,2$, and let $\cl_k$ be a 1-dimensional foliation    determined
by a line field transversal to $\cs_k$.  Then define the transversal 3-dimensional foliations
 $\cf_1 \equiv \cs_1\times \cl_2$ and $\cf_2\equiv \cl_1\times \cs_2$ on $M_1\times M_2$.
 By a theorem of J. Wood [Wo] any 2-plane field
on a 3-manifold is homotopic to a integrable one.   Any such foliation can then be modified by 
introducing many Reeb components.  The 1-dimensional foliations can often be constructed to have
dense orbits.

\medskip

\Remark { \E.3. ($\bbd$-submanifolds)}  A $\bbd$-submanifold is a submanifold whose
tangent spaces are $\bbt$-invariant  with $\pm1$-eigenspaces of the {\bf same} dimension.
This second condition is not automatic, but if it holds at a point, it holds on the 
connected component of that point.  These subspaces are automatically integrable,
with leaves given by intersection with the leaves of the ambient manifold.\medskip
The graphs of  a $\bbd$-holomorphic mapping $f:X\to Y$ is a $\bbd$-submanifold.
So also is $f^{-1}(p)$ if $p$ is a regular value of $f$.

We now address the question of $\bbd$-line bundles and $\bbd$-vector bundles.
A {\sl $\bbd$-line bundle} on an arbitrary manifold is a family of free $\bbd$-modules
over $X$ which is locally isomorphic to $U\times \bbd$ with transition functions which
are smooth $\bbd^+$-valued functions.

\Def{ \E.4}  A {\sl holomorphic line bundle} on a $\bbd$-manifold is a $\bbd$ line
bundle whose transition functions can be chosen to be $\bbd$-holomorphic.

\medskip

Let $\cd$ denote the sheaf of germs of $\bbd$-holomorphic functions on a $\bbd$-manifold
$X$, and let $\cd^+$ denote the  sheaf of germs of holomorphic  $\bbd^+$-valued functions
on $X$.  Then   the \v Cech cohomology group 
$$
H^1(X,\cd^+)
$$
represents the isomorphism classes of holomorphic $\bbd$-line bundles on $X$.
Note that the sheaf  sequence $0\to \cd \harr{\exp}{} \cd^+\to 1$ is exact, and
therefore $H^1(X,\cd) \cong H^1(X,\cd^+)$. Now in local null coordinates we have 
$\cd = e\cdot \ce_u + \bar e\cdot \ce_v \cong \ce_u \oplus \ce_v$ where
$\ce_u, \ce_v$ are the sheaves of germs of functions of $u$ and $v$ respectively. 
Thus $H^1(X,\cd^+) \cong H^1(X, \ce_u) \oplus  H^1(X, \ce_v)$.

\Ex{ \E.5. (The Canonical Bundle)}   The bundle $\kappa \equiv \L^{n,0}$ of holomorphic
$n$-forms has transition functions $\det_{\bbd}  \! \left({\partial z^\a\over \partial z^\b}\right)$
which are holomorphic with values in $\bbd^+$.  In local $z$-coordinates a holomorphic  $n$--form can be written as 
$\Phi= F(z)  d z= F(z) d z_1\wedge \cdots \wedge d z_n$ where $F$ is a $\bbd$-holomorphic function. 
In null coordinates
one has $$\Phi= e\,  f(u) du + \bar e\, g(v) dv.\eqno{( \E.4)}$$
From this we see that on a $\bbd$-manifold which is a product $X=M_1\times M_2$ as above, there always
exists a global holomorphic $n$-form which is nowhere-null, i.e., $\Phi\wedge \bar \Phi$ is never zero.
(For elements of $\L^{n,0}\dn$ the concepts of null and non-null make sense without a metric).

\Remark{\E.6. (Holomorphic $\bbd$-Vector Bundles)}
There is also the notion of  a holomorphic $\bbd$-vector bundle of higher rank $m$ where one has local
trivializations $U_\a\times \bbd^m$ and transition functions which are holomophic maps
to ${\rm GL}_m^+(\bbd) \equiv \{A\in\M_{m,m}(\bbd) : \det_\bbd(A) \in\bbd^+\}$.
Note that  one can twist the Dolbeault complex by such bundles.

  \vskip.3in

\centerline{\headfont \ \F.\  Hermitian D-Manifolds. }  
\medskip

Let $X$ be a $\bbd$-manifold and let $TX = N^+ \oplus N^-$ be the decomposition
into the tangent spaces of the two foliations, i.e, the $\pm 1$-eigenbundles of the  
endomorphism $\bbt: TX\to TX$  in (\E.3).

\Def{\F.1} By a {\sl hermitian metric} on $X$ we mean a non-degenerate $\bbd$-valued real-bilinear form
$(\cdot, \cdot)$ on the fibres of $TX$ such that
$$
(V,W)\ =\ \overline { (W,V)}
\and 
(\bbt V,W) \ =\ \tau (V,W) \ =\ -(V, \bbt W)
\eqno{(\F.1)}
$$
for  all $V,W\in T_xX$ at each point $x$.
\medskip

Note that 
$$
(\bbt V, \bbt W ) \ =\ -(V,W)
\eqno{(\F.2)}
$$
from which it follows that 
$$
N^+ \ \ {\rm and}\ \  N^- \ \ {\rm are\ null \ spaces \ for\ \ } (\cdot, \cdot).
$$
Recall (3.3) which expresses $(\cdot, \cdot)$ in terms of its  real and imaginary parts:
$$
(V,W) \ =\ \langle V,W \rangle - \tau \, \o(V,W) \qquad
$$
and note that $\langle V,W \rangle\ =\ \langle W,V \rangle$,  $\o(V,W)\ =\ -\o(W,V)$ and 
$$
 \langle V,  \bbt W \rangle\ =\ \o(V,W).
\eqno{(\F.3)}
$$
Thus $\langle V,W \rangle $ is a semi-riemannian metric of type $(n,n)$.  The spaces $N^\pm$ are null
for this metric and for its associated ``K\"ahler'' form $\o$. Both of these forms are characterized by the same 
bundle isomorphism:
$$
A: N^+\ \harr {\cong}{}\  (N^-)^*.
\eqno{(\F.4)}
$$
Writing $V=(V^+,V^-)$ with respect to the decomposition  $TX = N^+ \oplus N^-$, we have
$\langle V,W \rangle = A(V^+)(W^-) + A(W^+)(V^-)$ and 
$\o(V,W) = A(V^+)(W^-) - A(W^+)(V^-)$. In terms of local $z$-coordinates and local null coordinates 
$(u,v)$ one has
$$\eqalign
{
ds^2\ &=\ \sum_{j,k=1}^n a_{jk} dz_j\otimes d\bar z_k\ =\ \sum_{j,k=1}^n \tilde a_{jk} du_j\circ d v_k, \qquad{\rm and}\cr
\o\ &=\ \sum_{j,k=1}^n a_{jk} dz_j\wedge d\bar z_k\ =\ \sum_{j,k=1}^n \tilde a_{jk} du_j\wedge d v_k
}
$$
{\bf Note \F.2.} The isomorphism (\F.4) shows that the existence of a $\bbd$-hermitian metric on $X$ puts further topological restrictions on the bundle $TX$ (as opposed to the complex hermitian case).
For example, any product $X=M_+\times M_-$ of two $n$-manifolds is a $\bbd$-manifold, but
  (\F.4) shows that a hermitian metric exists on $X$ if and only if $M_+$ and $M_-$ are both  parallelizable.
  Also on any manifold $X$ the condition $TX\cong N\oplus N$ is restrictive.

\Note{\F.3}
A $\bbd$-submanifold $Y\ss X$ of a hermitian $\bbd$-manifold may not be hermitian
in the induced metric.  One needs the additional hypothesis that $Y$ is a symplectic
submanifold, i.e., that $\o\bigr|_Y$ is non-degenerate on $Y$.

\Note{\F.4} On any hermitian $\bbd$-manifold there exists the canonical Levi-Civita connection for the
semi-riemannian metric.  There also exists a canonical hermitian connection, characterized by the
fact that the metric and $\bbt$ are parallel and $\nabla^{0,1} = \dbar$.

\Remark{\F.5} 
In parallel with  the complex case, canonical
connections exist on {\sl any} holomorphic $\bbd$-bundle equipped with a hermitian metric
(cf. Remark \E.6).
The proof follows the complex case, where one computes the connection 1-form in a holomorphic
frame as $\o= h^{-1}\partial h$ where $h$ is the matrix determined by the frame and the metric.
The curvature 2-form $\Omega = d\o-\o\w\o$ transforms to $A\O A^{-1}$ under a GL$_n(\bbd)$-valued
frame change $A$, so that det$_\bbd(I+\tau \O)$ is a globally defined form  on $X$, called the
{\sl total $\bbd$-Chern form of the bundle}. More generally, for any other GL$_n(\bbd)$-invariant
polynomial $\phi$, the $\phi$-Chern form of the bundle, $\phi(\O)$, is a globally defined $d$-closed
differential form.  The standard transgression formula holds, and the cohomology class of $\phi(\O)$
in $H^*(X,\bbd)$ is independent of the hermitian metric (or connection) on the bundle.
Finally, if $h=e\cdot g+\bar e\cdot g^t$ defines a matrix $g$ with values in  GL$_n(\bbr)$, 
then one can show that in null coordinates 
$$
\tau \O = \tau \dbar(\partial h \cdot h^{-1}) = \bar e d_u \left (   (d_v g) \cdot  g^{-1}  \right) - 
e d_v \left (   (d_u g)  \cdot g^{-1}  \right)
$$

\Remark{\F.6} 
Suppose $X$ carries a nowhere-null holomorphic $n$-form  $\Phi$ (see Example  \E.4). Then changing the
metric by the conformal factor $\|\Phi\|^{1\over n}$ gives a new hermitian metric in which
$\Phi$ has constant length.  If $\Phi'$ is another holomorphic $n$-form with 
constant length in this new metric, then $\Phi'=\a \Phi$ for $\a\in\bbd$ with $\a\bar\a\neq0$.

  \vskip .3in
 
 \def\Kahler{K\"ahler }

\centerline{\headfont \ \G. \  K\"ahler D-Manifolds. }  
\medskip

We now consider the following natural class of hermitian double manifolds.

\Def{\G.1} A hermitian $\bbd$-manifold is said to be {\sl K\"ahler } if  the form $\o$ is $d$-closed.

\medskip
Interestingly, all the standard characterizations of complex K\"ahler manifolds carry over to this context.
For example, we have the following.

\Prop{\G.2} {\sl
A hermitian $\bbd$-manifold $X$ is  K\"ahler
if and only if the canonical hermitian connection on $TX$ agrees with the Levi-Civita connection of the 
semi-riemannian metric
}
\Prop{\G.3} {\sl   Let $X$   be a hermitian $\bbd$-manifold with K\"ahler form $\o$ and
structure map $\bbt$. Then $X$ is K\'ahler  if and only if $\bbt$ is parallel in the Levi-Civita connection.
}\medskip

In the  last assertion one does not need to 
assume {\sl apriori} that the subbundles $N^\pm$  are integrable. The hypothesis $\nabla \bbt=0$ implies
that $N^\pm$  are integrable. 

\Note{\G.4} In the null-coordinate approach,  a K\"ahler $\bbd$-manifold is simply a symplectic manifold
with a pair of transverse Lagrangian foliations.  For this reason  K\"ahler $\bbd$-manifolds are sometimes
referred to in the literature as {\sl bi-Lagrangian manifolds}.

\Note{\G.5} Any symplectic $\bbd$-submanifold of a  \Kahler $\bbd$-manifold, is \Kahler
in the induced metric.

\medskip
\noindent
{\bf Example \G.6.(Surfaces).}  When  $n=1$ a \Kahler $\bbd$-manifold is simply a
 surface equipped with a pair of 1-forms
$\a$ and $\b$ such that $\a\wedge\b$ is nowhere vanishing.  Here the metric is $ds^2=\a\circ \b$ and the \Kahler form is $\o=\a\wedge\b$.
Equivalently it is a Lorentzian surface.

  \vskip .3in

\centerline{\headfont \ \H. \  Ricci-Flat K\"ahler  D-Manifolds. }  
\medskip

  Let $X$ be a hermitian $\bbd$-manifold and suppose $\Phi$ is a nowhere-null holomorphic section of the
  canonical bundle $\kappa = \L^{n,0}(X)$.  Then the real 2-form
  $$
\O\ \equiv\   \tau \partial\dbar \log \|\Phi\|  
 \eqno{(\H.1)}
  $$
  is independent of $\Phi$.
To see this note first that in local holomorphic coordinates $z$ we have $\Phi = a(z) dz_1\wedge\dots\wedge dz_n$,
where $a(z)$ is a $\bbd$-holomorphic function, and $\|\Phi \|^2 =  \lambda(z) a(z)\bar a(z)$ where $\l>0$.
Hence we have   $\tau  \partial\dbar \log \|\Phi\|^2  = \tau  \partial\dbar \log \lambda$. Changing coordinates
to $z'$ gives a new metric coefficient $\l' = \l F\overline F$ where $F$ is a holomorphic function
(the determinant of $\partial z/\partial z'$). Thus $ \partial\dbar \log \lambda =  \partial\dbar \log \lambda'$.
This form $\O$ is the curvature of the canonical hermitian connection on $\kappa$.  (This generalizes to any hermitian holomorphic line
bundle.) 

We observed in Remark \F.6 that given $X$ and $\Phi$ as above, there is a conformally
equivalent metric on $X$ for which $\|\Phi\|\equiv 1$. Evidently, in this metric the curvature
$\O\equiv 0$.  

As remarked in (\F.4), if one passes to  null coordinates, then 
$\Phi\ =e\, f(u) du+\bar e \,g(v)dv.$ For this reason $\bbd$ geometry plays
an important role in the mass-transport problem.

When $X$ is  \Kahler we have the following.

 \Prop{\H.1} {\sl
 Let $X$ be a   K\"ahler $\bbd$-manifold, and $\O$ the curvature 2-form of the canonical
 line bundle as above.  Then
 $$
 \O(V,W)\ =\ -Ric(V, \bbt W)
 \eqno{(\H.2)}
 $$
 where $Ric$ is the Ricci curvature tensor of the Levi-Civita connection on $TX$.
 In particular,  the canonical bundle is flat if and only if
 the Ricci curvature of $X$ is zero.
 } 
\pf   The argument goes precisely as in the complex case, with $i$ replaced by $\tau$ and
the almost complex structure $J$ replaced by the structure map $\bbt$.\medskip

Thus a Ricci-flat   K\"ahler $\bbd$-manifold is the exact analogue of a  Calabi-Yau manifold.
The holonomy will lie in the group ${\rm SU}_n(\bbd)$ (cf. [B].)

These manifolds can be characterized, in analogy with Hitchin's [Hi$_1$] 
description of the complex case, as follows.

\Prop{\H.2} {\sl
A Ricci-flat   K\"ahler $\bbd$-manifold is equivalent to the data of a
symplectic $2n$-dimensional manifold $(X,\o)$ together with two 
$d$-closed real $n$-forms $\phi$, $\psi$ such that:
\medskip

(1) \ \  $\Phi = \phi+\tau\, \psi$ is a simple (indecomposable)  $\bbd$-valued $n$-form.

\medskip

(2) \ \  $\Phi \wedge \o =0$, i.e. $\phi\wedge \o =\psi\wedge \o =0$  
\medskip

(3) \ \  
$$
\Phi\wedge \overline \Phi \ =\ \cases
{
\quad\  \o^n \qquad {\sl for \ \ } n \ \  {\sl even}  \cr
-\tau \o^n\qquad {\sl  for \ \ } n \ \ {\sl odd}.  \cr
}
$$
}
This proposition can be restated as follows.

\Prop{\H.2$'$} {\sl
A Ricci-flat   K\"ahler $\bbd$-manifold is equivalent to the data of a
symplectic $2n$-dimensional manifold $(X,\o)$ together with two 
$d$-closed real $n$-forms $\a$, $\b$ such that:
\medskip

(1) \ \  $\Phi = e\,\a+\bar e\,\b$ and the real $n$-forms $\a$ and $\b$ are simple.

\medskip

(2) \ \  $\a\wedge \o =\b\wedge \o =0$  
\medskip

(3) \ \  
$
\a\wedge\b\ =\ \o^n$.
}

\medskip
To see the equivalence of these statements, set $\a =  \phi-\psi$ and $\b= \phi+\psi$
and note that $(\a_1\wee \a_n)\cdot e  + (\b_1\wee \b_n )\cdot \bar e  = 
(\a_1\cdot e +\b_1\cdot \bar e)\wee  (\a_n\cdot e +\b_n\cdot \bar e)$

\def\psii{\x}

There are several interesting examples.

\Ex{\H.3}  By Proposition 11.2$'$  a Ricci-flat  \Kahler $\bbd$-manifold of (real) dimension 2 is  simply a   surface
$\Sigma$ equipped with a pair of {\sl closed} 1-forms $\a$, $\b$ such that $\a\wedge\b$ 
never vanishes.  (The K\"ahler form is $\o=\a\w\b$.)

Note that given this data, $\psii\equiv \a+i\b$ defines a conformal structure on $\Sigma$
in which $\psii$ is a (complex) holomorphic 1-form. Conversely if $\Sigma$ is a Riemann
surface with a nowhere vanishing holomorphic 1-form $\psii = \a+i\b$, then $(\Sigma,\a,\b)$ is
a Ricci-flat \Kahler $\bbd$-manifold.  This leads to the following.

\Prop{\H.4} {\sl
Every non-compact connected surface carries a Ricci-flat K\"ahler $\bbd$-structure.
The only compact connected surface to carry such a structure is the torus.
}
\pf
Any surface $\Sigma$ can be given a conformal structure. If $\Sigma$ is a  non-compact connected Riemann surface, then it  carries a  nowhere vanishing  holomorphic 1-form  since it is Stein. 
If $\Sigma$ is compact, it  carries such a  1-form if and only if its genus is zero.  \qed

\medskip

Since products of Ricci-flat  \Kahler $\bbd$-manifolds are again of this type, we obtain
examples with large fundamental groups in all dimensions.

\def\a{a}\def\b{b}

 \Ex{\H.5}  Let $\Sigma, \Sigma'$ be  Riemann 
surfaces endowed with nowhere-vanishing holomorphic 1-forms
$$
\psii\ =\ \a+i\b \and \psii'\ =\ \a'+i\b'.
$$
 On $X=\Sigma\times \Sigma'$ we have the symplectic form and metric
given by
$$
\o \ =\ \a\wedge \a' + \b\wedge \b'
\and
ds^2 \ =\ \a\circ \a' + \b\circ \b'.
$$
The tensor $\bbt:T(\Sigma\times\Sigma') \to T(\Sigma\times\Sigma')$ is defined by 
$
\bbt\ =\ \left(\matrix{I&0\cr 0&-I}\right).
$
The holomorphic 1-form
$$
\Phi \ = \  e\cdot \a\wedge \b + \bar e\cdot \a'\wedge \b' \ \equiv \  e\cdot \alpha  + \bar e\cdot\beta
$$
satisfies conditions (1), (2) and (3) above, and so $X$ is  a Ricci-flat \Kahler $\bbd$-manifold.
 One can also take products of such 4-manifolds.

\def\a{\alpha}\def\b{\beta}

 \Remark{\H.6. (The Calabi Question)} Suppose $(X,\o)$ is a K\"ahler $\bbd$-manifold with trivial canonical bundle.
It  is natural to ask the {\sl Calabi Question}:  {\sl Does there exist a  K\"ahler form $\o'$ on the $\bbd$-manifold
$X$, which is cohomologous to $\o$ and is Ricci-flat?}\
 Interestingly, the answer is: {\sl Not always}.  Consider the $\bbd$-structure on the torus $S^1\times S^1$
 given by the transverse foliations $\cf_1$ and $\cf_2$, where $\cf_1$ consists of the  circles 
 $S^1\times \{\theta\}$ for $\theta\in S^1$, and $\cf_2$ has exactly one compact leaf $L=\{\theta_0\}\times  S^1$
 (whose complement is foliated by leaves which spiral barber-pole fashion from one side of $L$ to the other). Any choice of nowhere vanishing 1-forms $\a$ and $\b$ 
so that $\a$ vanishes on $\cf_1$ and $\b$ vanishes on $\cf_2$ makes this a K\"ahler $\bbd$-manifold.
However, {\sl it is not possible to choose $\a$ and $\b$  to be closed, that is, A Ricci-flat K\"ahler structure
does not exist.} To see this, suppose $\b$ were closed, and consider a thin strip $S\equiv (\theta_0 -\e,\theta_0+\e)\times S^1$ about the compact leaf $L$. Since $\b\bigr|_{L}=0$, we conclude that $\b\bigr|_{S}$  is exact
by the deRham Theorem.  That is, $\b=df$ where $f:S\to \bbr$ is a function whose level sets are the leaves
of $\cf_2$ in $S$. However, since $df=\b\neq 0$
on $S$, one sees that the level sets of $f$ near $L$ must be compact (circles).

Taking products gives examples in all dimensions.

 \vskip.3in

\centerline{\headfont  12.  Split SLAG Submanifolds in the General Setting\ }  
\medskip

Let $X$ be a Hermitian $\bbd$-manifold of $\bbd$-dimension $n$ with a nowhere null holomorphic
$n$-form  $\Phi$. Let $\o$ be the K\"ahler form of the hermitian metric. (We do not require $d\o=0$.)

\Def{12.1}  
An real oriented $C^1$ submanifold $M$ of dimension $n$ in $X$ is said to be 
{\bf split SLAG}  (or $\wt {\bf SLAG}$) if
\medskip

(1)\ \ $M$ is Lagrangian, i.e., $\o\bigr|_M=0$.

\medskip

(2)\ \   ${\rm Im}\, \Phi \bigr|_M=0$.

\medskip

(3)\ \ $M$ is positive space-like.

\Theorem{12.2} {\sl
Suppose $\|\Phi\|\equiv 1$ on $X$.  Then any compact split SLAG submanifold with boundary
$(M,\partial M)$ in $X$ is homologically volume-maximizing, i.e., 
$$
{\rm vol}(M)\ \geq\ {\rm vol}(N)
$$
for any other positive space-like submanifold $N$ such that  $\partial N = \partial M$
and $M-N$ is homologous to zero in $X$.

Equality holds if and only if $N$ is also split SLAG.
}
\pf The same as the proof of Theorem  \SS.3.\qed

\Remark{12.3}  One can always change the metric on $X$ by a conformal factor so that
$\|\Phi\| \equiv 1$.

\Remark{12.4}  Theorem 12.2 carries over to split SLAG subvarieties.  See Appendix D.

 \vskip.3in

\centerline{\headfont  \QQ.  Deformations and the McLean Theorem in the Split Case.} \medskip

In his   Duke University thesis [Mc$_{1,2}$], R. McLean proved that the moduli space of Special Lagrangian
submanifolds of a Ricci-flat K\"ahler manifold is smooth, and its tangent space at a point (submanifold)
$M$ is canonically identified with the space ${\bf H}^1\cong H^1(M;\bbr)$ of harmonic 
1-forms on $M$. This result carries over to the split Special Lagrangian case (see Warren[W$_2$]). 
We present this result
and sketch the  proof based in part on the notes of Hitchin [Hi$_1$].

\Theorem {\QQ.1}  {\sl Let $X, \o, \Phi$ be as in \S 12 and assume $d\o=0$ (so $X$ is a Ricci-flat 
K\"ahler $\bbd$-manifold). Let $M\ss X$ be a split SLAG submanifold, and suppose $M_t, -\e<t<\e$, is a variation
 through split SLAG submanifolds with $M=M_0$.  Then the corresponding normal variational vector field
$\nu$ canonically gives a   harmonic 1-form $\theta$ on $M$ by setting
$$
\theta\ =\   \nu\hk \o.
$$
If $M$ is compact, then,  in a neighborhood of $M$, the moduli space $\gerM$ of split SLAG submanifolds
of $X$ is a manifold of dimension $b_1(M)$ whose tangent space at  $M' \in \gerM$ is canonically
identified with the space of harmonic {1-forms} on $M'$.  In particular, this endows $\gerM$ with
a natural riemannian metric.
}

\pf
Let $F:M\times (-\e, \e) \to X$ be a smooth variation of $M=M_0$ so that $M_t \equiv F(M\times \{t\})$
is positive space-like for all $t$. Then by Definition 12.1, this is a variation through  split SLAG submanifolds if and only if
$$
F^*\o\bigr|_{M\times \{t\}} \ =\ 0 \and
F^*({\rm Im}\Phi)\bigr|_{M\times \{t\}} \ =\ 0  \fa t,
$$
or equivalently
$$
F^*\o = \wt\theta \wedge dt 
\and
F^*({\rm Im}\Phi) = \wt\vf \wedge dt 
$$
where $\wt \theta$ and $\wt \vf$ are a 1-form and an $(n-1)$-form respectively on $M\times (-\e, \e) $.
Note that the restrictions 
$\theta = \wt\theta\bigr|_{M\times \{t\}}$ and $\vf = \wt\vf\bigr|_{M\times \{t\}}$
are independent of the choice of $\wt \theta$ and $\wt\vf$.
Since $dF^*\o = dF^*({\rm Im}\Phi) = 0$,  these restrictions satisfy $d\theta  = d\vf = 0$.
Note also that $\theta$ and $\vf$ correspond to the restriction to $M_t$ of the contraction
of $\nu \equiv F_*({\partial \over \partial t})$ with $\o$ and Im$\Phi$ respectively. In fact, by 12.1 (1) and (2), only
the normal part of $F_*({\partial \over \partial t})$ survives under restriction.
The first part of the Theorem is now a consequence of the following algebraic assertion:
$$
\vf\ =\ *\theta
\eqno{(\QQ.1)}
$$
where $*$ denotes the Hodge star-operator  with respect to the induced riemannian
metric on $M_t$. The proof of  (\QQ.1)  is straightforward and  left to the reader.

The proof of the unobstructedness of the moduli space near $M$, when $M$ is compact, 
follows exactly the argument given in [Mc$_2$].\qed

 \vskip.3in

\centerline{\headfont  \MM.  Relation to the Mass Transport Problem --} \medskip
\centerline{\headfont Work of Kim, McCann and Warren\ }  
\medskip

We first recall the classical Monge-Kantorovich Optimal Mass Transport Problem
in a smooth setting. Fix $n$-dimensional manifolds $U$ and $V$ with domains
$\O_U\ss\ss U$, $\O_V\ss\ss V$ and smooth positive densities $\rho_U, \rho_V$
on $\O_U$ and $\O_V$ respectively.  Let $c:\O_U\times \O_V\to \bbr$ be a smooth
``cost'' function. Set
$$
\cm \ \equiv\ \{F:\O_U\ \to\ \O_V : F_*(\rho_U) =\rho_V\}
$$

\medskip\noindent
{\bf Problem.}  Find a mapping in $\cm$ which minimizes the functional
$$
C(F)\ \equiv\ \int_{\O_U} c(u,F(u))\, \rho_U
$$
 It turns out that this is equivalent to the following.
 Consider the space
 $$
 \cl\ \equiv\  \{(f,g) \in C(\O_U)\times C(\O_V) : f(u)+g(v)\leq c(u,v) \ \forall\ u,v  \}
 $$
 
\medskip\noindent
{\bf Dual Problem.} Find a pair $(f,g)\in \cf$ which maximizes the functional
$$
J(f,g)\ \equiv\ \int f(u)\, \rho_U \ +\ \int g(v) \, \rho_V \
$$

Consider now the case $U=V=\rn$. 
Under appropriate conditions on the cost function these equivalent problems have a 
 unique solution.  For the classical cost function $c(u,v) =\half |u-v|^2$ the solutions are related by 
 $$
 F\ =\  df\ =\ (dg)^{-1} \ =\ F^{-1}.
 $$
 Note that $f$ and $g$ are Legendre transforms of one-another. 
If we   write $\rho_U = \rho(u) du$ and $\rho_V = \wt{\rho}(v) dv$ for smooth positive functions
 $\rho, \wt \rho$, 
 then $f$ satisfies the Monge-Amp\`ere equation
 $$
 \det \, \Hess f\ =\ { \rho \over \wt \rho(df)}.
 $$
 For a complete and detailed exposition of the modern state of knowledge
 pertaining to the optimal transport problem, the reader is referred to [V].
 
 Even when one is working in the general case presented above, there are 
 existence and uniqueness results.     Kim, McCann and Warren have found
 that this problem and its solution fit beautifully into the framework 
 of K\"ahler $\bbd$-manifolds.  Briefly it goes as follows.
 
 Consider the $\bbd$-manifold
 $$
 U\times V.
 $$ 
The first  assumption is  that the given cost function gives a {\bf global K\"ahler potential},
 i.e., 
 $$
 \o \ \equiv\ d_u d_v c
 $$
 is a symplectic (i.e., non-degenerate) 2-form on $U\times V$.
 (Otherwise said, the matrix ${\partial^2 c\over\partial u_i\partial v_j}$ 
 is everywhere non-singular,)
 This makes $U\times V$ into a {\bf  K\"ahler $\bbd$-manifold}.
 One also imposes a so-call ``twist condition''  on $c$ which says that
 the maps $u\mapsto (d_v c)(u,v)$ and $v\mapsto (d_u c)(u,v)$
 are injective, and which we can ignore here.
 
 Observe now that our K\"ahler $\bbd$-manifold is equipped with a {\bf holomorphic
 $n$-form}
 $$
 \Phi \ \equiv \ \rho(u)\, du\cdot e \ +\ \wt\rho(v)\, dv\cdot\bar e
 $$
 \Theorem{\MM.1. [KMW]}  {\sl
 The graph of the unique solution to the Monge-Kantorovich optimal mass transport problem
 given above is a split SLAG submanifold of the K\"ahler $\bbd$-manifold $U\times V$.
 }
 \medskip
 
 Note that in general $U\times V$ is not Ricci-flat.  Hence the graph of $F$ is only volume-maximizing
 for a conformally related metric (cf. Remark 12.3).

  \vskip .3in

\centerline{\headfont  \TT. Further Work of Mealy} \medskip
 In this section we give a brief overview of the primary work on semi-riemannian calibrated 
 geometry [M$_{1,2}$].  Let $\bbr^{p,q}$ denote the standard semi-euclidean vector space of signature $p,q \geq 1$.
 Choose a connected component $G$ of the grassmannian of oriented non-degenerate  linear subspaces of 
 $\bbr^{p,q}$ of dimension $k=r+s$ and signature $r,s$.  There is a canonical embedding of $G$
 into $\L^{r+s}\bbr^{p,q}$ via the wedge product of any oriented orthonormal basis for the subspace.
 
 \Def{\TT.1} Given a subset $A\ss G$, a form $\phi\in \L^{r+s}( \bbr^{p,q})^*$ is a {\bf calibration of type $A$} 
 if
 $$
 \phi(\x)\ \geq\ 1 \fa \x\in A.
 \eqno{(\TT.1)}
 $$
 Moreover, the set 
 $$
 G(\phi) \ =\ \{\x\in A: \phi(\x)=1\}
 $$
of calibrated planes is usually assumed to be non-empty.
 An oriented submanifold $M$ is of {\bf type $A$} if each 
 tangent plane is in $A$.  We say $M$ is {\bf calibrated} by $\phi$ if $M$ is of type $G(\phi)$.
 \medskip
 
 The cases where $A$ can be taken to be a full connected component of $G$ are very limited.
 Ignoring the two cases of dimension $k=0$ and dimension $k=n$, where $\phi$ is $\pm$ the volume form
 on $\bbr^{p,q}$, and eliminating duplications from interchanging $p$ and $q$, Mealy shows that
 only two cases (essentially one) remain. The Grassmannian $G_{\rm space}(p, \bbr^{p,q})$ of maximally space-like oriented planes in $\bbr^{p,q}$, has two connected components.
 Let $G^+_{\rm space}(p, \bbr^{p,q})$ denote the component containing $\oa{\bbr}^p = \bbr^p\times \{0\}$
 with its standard orientation.  The other component is 
 $G^-_{\rm space}(p, \bbr^{p,q}) = -G^+_{\rm space}(p, \bbr^{p,q})$.  Henceforth we only consider
 $G^+_{\rm space}(p, \bbr^{p,q})$.  Thus the problem of finding semi-riemannian calibrations $\phi$ of type
 $A$, where $A$ is a connected component of the non-degenerate Grassmannian, is reduced to the case
 $A\equiv G^+_{\rm space}(p, \bbr^{p,q})$ and the problem of finding a form
 $$
 \phi\in \L^{r+s}( \bbr^{p,q})^*\qquad{\rm with}\ \ \phi(\x) \geq1\quad \forall\, \x\in \phi\in  G^+_{\rm space}(p, \bbr^{p,q}).
  \eqno{(\TT.2)}
 $$

 \Remark{\TT.2} In the next section we will give an important example of ``type $A$'' 
 where $A$ is a proper subset of a component of the Grassmannian.
 
 \Theorem{\TT.3}  {\sl
 Suppose that $\phi$ is a calibration of type $A$ and $M$ is a compact  submanifold of $\bbr^{p,q}$
 calibrated by $\phi$. If $N$ is any other  compact submanifold of type $A$ with $\partial N=\partial M$, 
 then
 $$
 {\rm vol}(M)\ \geq \ {\rm vol}(N)
 $$
 with equality if and only if $N$ is also calibrated by $\phi$.
 }
 \pf See the proof of Theorem \SS.3.\qed

 \Ex{1. (Point Calibrations)}  The inner product $\bra \cdot \cdot$ of signature $p,q$ on $\bbr^{p,q}$ induces an isomorphism 
 $\bbr^{p,q} \harr \cong \ (\bbr^{p,q})^*$ which extends to an equivalence 
 between $\L^k \bbr^{p,q}$ and $\L^k ( \bbr^{p,q})^*$.  Given $\xi\in G_{\rm space}^+(p,  \bbr^{p,q})$,
 the equivalent form $\phi_\x\in \L^p( \bbr^{p,q})^*$ defined by $\phi_\x(\eta)=\bra\x\eta$
  is a calibration of type
 $G_{\rm space}^+(p,  \bbr^{p,q})$.  More precisely one has:
$$
\eqalign
{
\phi_\x(\eta)\ \geq\ 1 \ \  {\rm for}\ \  & \eta \in G_{\rm space}^+(p,  \bbr^{p,q})   \cr
{\rm with\  equality }\ \ &\iff\ \ \eta\ =\ \xi
}
 \eqno{(\TT.3)}
 $$
 
 With a conflict of terminology, the case $p=1$ yields a calibration proof of the twin paradox for curves in 
 $\bbr^{1,q}$ by taking the first axis to be time (see [H, Prop 4.19]).  
 The proof of the inequality above  is immediate (see the bottom of page 797 in [M$_2$]) from a
 canonical form for $\eta$ with respect to $\xi$ under the action of the orthogonal group O$(p,q)$
 given on page 796 of [M$_2$].

 \Ex{2. (Divided Powers of the K\"ahler Form)} 
  Now take $\bbc^{p,q}$ to be $\bbc^{p+q}$ equipped with the standard complex hermitian  inner
 product of signature $p,q$, and let $\o$ denote the corresponding ``K\"ahler'' form.  In general,
 $\phi \equiv {1\over k!}\o^k$ is {\sl not} a calibration.  It is a calibration only when $k=p$.

 \Ex{3. (Split Associative, Coassociative, etc.)} 
Replacing the octonions $\bbo$ by their split companion $\wt\bbo$ (analogous to replacing
 $\bbc$ by its split companion $\wt \bbc = \bbd$), the associative, coassociative and Cayley calibrations
 have counterparts that are calibrations. 
 \medskip
 
 Finally it is worth noting that all possible calibrations of type $G_{\rm space}^+(3,\bbr^{3,3})$
 are classified in [M$_1$]. This is the lowest dimensional non-trivial case.

  \vskip .3in

\centerline{\headfont  \UU.  Lagrangian Submanifolds of Constant Phase}\smallskip
\centerline{\headfont  and Volume Maximization}
 \medskip
 
 In this section we look at the general picture of all oriented Lagrangian submanifolds with non-degenerate
 metric.  We show that such  submanifolds are minimal (mean curvature zero)  if and only if they are of 
 {\sl constant phase}.  When and only when the submanifolds are also purely space-like or time-like,
they are homologically volume-maximizing among purely space-like or time-like submanifolds 
of the same type.  However, in all other signatures they are homologically  volume-maximizing
 among Lagrangian submanifolds of the same type.
For convenience we shall work in $\dn$, however the  results hold in any Ricci-flat $\bbd$-manifold.

The first result of this nature was noted in [HL$_1$, page 96] in the standard special Lagrangian setting.
The phase $\theta$ of a Lagrangian submanifold $M$ of $\bbc^n$ was defined using the fact that 
$dz(\oa M)=e^{i\theta}$ is of unit length (Proposition 1.14 and (2.18)).  For any tangent vector $V$
to $M$ one has that $V(\theta) = \bra {JV} H$  (equation (2.19))
where $H$ is the mean curvature of $M$.

Dong [D] considered the case of non-degenerate Lagrangian submanifolds $M$ of $\bbc^{p,q}$ ($\bbc^{p+q}$ 
with the standard hermitian inner product of signature $p,q$).  Again the phase $\theta$ can be defined
by $e^{i\theta} = (dz)(\oa M)$ since it is of unit length. Despite the fact that Re$\,dz$ is not a calibration,
Dong noted that $V(\theta) = \bra {JV}H$ still holds, so that $M$ is of mean curvature zero
if and only if the phase function is constant ([D, Lemma 2.1]).
Each non-degenerate Lagrangian subspace in $\bbc^{p,q}$ must have signature $p,q$.
We would like to make the further observation  that  Re$\,dz$, although  
not a calibration in the pure sense, is a calibration of type $A=\LAG$, the
set of oriented non-degenerate Lagrangian subspaces, because $(dz)(\x) = e^{i\theta}(\x)$
for $\x\in\LAG$.   By arguing as in the proof of Theorem
\TT.4, this proves that a constant phase Lagrangian submanifold is volume-minimizing
among all other Lagrangian submanifold with the same boundary.

We now turn attention back to $\dn$.
Let $\LAG$ denote the set of oriented non-degenerate  Lagrangian $n$-planes in $\dn$.  
The set $\LAG$ decomposes into $2n+2$ connected components
$$
\LAG\ = \ \coprod_{p+q=n} \LAG_{p,q}^\pm
$$
where $\LAG_{p,q}^\pm$ consists of planes for which the induced metric has signature $(p,q)$ and orientation
$+$ or $-$ when compared to an {\sl a priori} chosen model. Each $\LAG_{p,q}^+$ and $ \LAG_{p,q}^-$
is an orbit of the group U$_n^+(\bbd)$, and the pair $\LAG_{p,q}^+ \cup \LAG_{p,q}^-$ is an orbit
of U$_n(\bbd)$.

\Note{\UU.1} The oriented model planes can be chosen so that 
$$
\pm \tau^q dz \bigr|_P\ =\  e^{\tau \theta_P} d{\rm vol}_P
 \qquad {\rm for} \ \  P\in  \LAG_{p,q}^\pm.
\eqno{(\UU.1)}
$$
where $\theta_P\in \bbr$ and $d{\rm vol}_P$ is the unit (positive) volume form on $P$.

\medskip
If $M$ is an oriented connected non-degenerate Lagrangian submanifold 
of signature $p,q$, then all  its tangent planes lie either in $\LAG_{p,q}^+$  or in $ \LAG_{p,q}^-$
depending on the orientation of $M$ and we say that $M$ is of {\bf type} $\LAG_{p,q}^+$
or of  {\bf type} $\LAG_{p,q}^-$. Therefore by (\UU.1)
we have that
$$
\pm \tau^q dz \bigr|_M\ =\  e^{\tau \theta} d{\rm vol}_M.
\eqno{(\UU.2)}
$$
\Def{\UU.2}  This  smooth real-valued function $\theta$ on $M$ 
will be called the {\bf phase function} on $M$.
 If $\theta$ is constant, we say that $M$ has {\bf  constant phase}.

\Prop {\UU.3}  {\sl 
Let $M\ss\dn$ be an oriented non-degenerate Lagrangian submanifold which is connected.  Then 
$M$ is a minimal (mean curvature zero) submanifold if and only if $M$ has constant phase.
}

\pf This follows immediately from the fact that for any tangent vector field $V$ on $M$
$$
V\cdot \theta \ =\ \bra {\bbt V} H
\eqno{(\UU.3)}
$$
where $H$ is the mean curvature vector of $M$. By (\F.3) this is equivalent to the statement that
$$
d\theta\ =\ H\hk \o
$$
A proof of (\UU.3) is given at the end
of the section.\qed

\Prop{\UU.4} {\sl The form
$$
\phi_{\theta} \ =\ \pm {\rm Re} \left\{ e^{-\tau\theta} \tau^q dz\right\}
\eqno{(\UU.4)}
$$
is a calibration of type $\LAG^\pm_{p,q}$.}

\pf Apply (\UU.1) and (\UU.2).\qed
\medskip

As noted in Theorem 5.3 and Remark 5.4, if $M$ is a compact submanifold (with boundary)
in $\dn$ of type $\LAG^+_{n,0}$ which has constant phase $\theta$, then for any submanifold
$N$ of type $\Gsp$ with the same boundary as $M$ we have vol$(M) \geq {\rm vol}(N)$, and
equality holds if and only if $N$ is also of type $\LAG^+_{n,0}$ with the same constant phase.
A similar result holds for $\LAG^+_{0,n}$ (and of course for $\LAG^-_{n,0}$),
but not for $\LAG^+_{p,q}$ if $p,q\geq 1$ unless one restricts attention of a subset
of the appropriate Grassmannian (See Remark \TT.2).

\Theorem{\UU.5} {\sl
Suppose that $(M,\partial M)$ is a compact submanifold  of type  $\LAG^\pm_{p,q}$ in $\dn$
which is of constant phase $\theta$.  If $(N, \partial N)$ is any other compact submanifold 
of the same type $\LAG^\pm_{p,q}$ with $\partial N=\partial M$, then
$$
{\rm vol}(M)\ \geq\ {\rm vol}(N)
$$
with equality if and only if $N$ is also of constant phase $\theta$.
}

\pf Apply Theorem \TT.3 and Proposition \UU.4.\qed

\Remark{\UU.6} This theorem carries over to any Ricci-flat K\"ahler $\bbd$-manifold.  
Here $M$ is only homologically volume maximizing, i.e., maximizing among $M'$ 
where $M-M'=\partial Y$ for some $(n+1)$-chain $Y$.

\medskip
\noindent
{\bf Proof of formula (\UU.3)}
Pick  a local  oriented orthonormal frame field $e_1,...,e_n$ on $M$.   Then
$$
dz(e_1\wee e_n)  \ =\ e^{ \tau  \theta}
$$
defines the smooth phase function $\theta$ for $M$.
Obviously  we have that $V(\theta) = \tau e^{-\tau\theta} V(e^{\tau\theta})
=   \tau e^{-\tau\theta} V(dz(e_1\wee e_n))$.
Now with $\nabla$ the riemannian connection on $M$, we have
$$
\eqalign
{
V\left( dz(e_1\wee e_n)  \right)  &=    
  dz \left ( \sum_{k=1}^n e_1\wee \nabla_V e_k \wee e_n\right)   \cr
&= dz \left ( \sum_{k=1}^n e_1\wee  \bra {\nabla_V e_k}{\bbt e_k} \bbt e_k \wee e_n\right)  
\ =\  \sum_{k=1}^n \bra {\nabla_V e_k}{\bbt e_k}  \tau e^{\tau\theta}\cr
}
$$
first because $\nabla dz=0$, second because $\bra {\nabla e_k}{e_k}=0$
and $dz(e_\ell \wedge \bbt  e_\ell \wedge \eta)=0$, with the third equality following
because $dz(e_1\wee \bbt e_k\wee e_n) = \tau dz(e_1\wee e_n) = \tau e^{\tau\theta}$.
Combining gives
$$
V(\theta)\ =\ \sum_{k=1}^n \bra {\nabla_V e_k}{\bbt e_k},
$$
but $\bra {\nabla_V e_k}{\bbt e_k} = \bra {\nabla_{e_k} V}{\bbt e_k} = -\bra {\nabla_{e_k} \bbt V}{ e_k}
= \bra {\bbt V}{\nabla_{e_k}e_k}$ and $H=\sum_k\left(\nabla_{e_k}e_k\right)^{\rm normal}$ is the mean
curvature.\qed

  \vskip .3in

\centerline{\headfont Appendix A }  
\medskip

\centerline{\headfont \ A Canonical Form for Space-like n-Planes in D$^n$.}  
\medskip

In this appendix we give a canonical form for space-like $n$-planes under the
unitary group.
We first treat the case $n=2$ which is the key.

\Lemma{A.1} {\sl
Given a space-like 2-plane $P$ in $\bbd^2$, there exist a space-like $\bbd$-unitary basis
$e_1,e_2$ for $\bbd^2$ and an angle $\theta\geq 0$ such that
$e_1$ and $\cosh \theta e_2 +\sinh\theta \bbt e_1$ form an orthonormal  basis for $P$.
}
\pf
  Pick any unit vector $e_1\in P$. Since $P$ is space-like, there
exists a unit time-like vector $\bbt e_2$ orthogonal to both $P$ and $\bbt e_1$.
Let this define the unit space-like vector $e_2$.  Then the unit vector $\e\in P$ orthogonal to $e_1$ belongs to the span of $e_2$ and $\bbt e_1$. Hence $\e= \cosh \theta e_2 + \sinh \theta \bbt e_1$, after possible sign 
changes of $e_1$ and $e_2$ to ensure $\theta\geq0$. \qed

\medskip
\noindent
{\bf Note.}  Using the  basis $e_1,e_2$ to define coordinates $z_1= x_1e_1+y_1\bbt e_1$ and
$z_2 = x_2e_2+y_2\bbt e_2$.
The  2-plane $P$ is defined by $y_1=  (\tanh\,\theta) x_2$ and $y_2=0$ (where $-1<\tanh\theta<1)$, i.e.,
$P$ is the graph over $\bbr^2$ in $\bbr^2\oplus \tau \bbr^2$ of the nilpotent matrix $A=\left(\matrix{0&\tanh\theta\cr
0&0}\right).$

\Prop{A.2. (Canonical Form)}  {\sl
Each space-like $n$-plane $P$ in $\bbd^n$ has an orthonormal basis 
$$
\{e_1, \ \cosh \theta_1 e_2 +\sinh\theta_1 \bbt e_1, \ ...\ ,
e_{n-1}, \ \cosh \theta_{n\over 2} e_n +\sinh\theta_{n\over 2}  \bbt e_{n-1}
\}  \qquad {\rm if\ } n {\rm \ is\ even, \ and}
$$
$$
\{e_1, \ \cosh \theta_1 e_2 +\sinh\theta_1 \bbt e_1, \ ...\ ,
e_{n-2}, \ \cosh \theta_{n-1\over 2} e_{n-1} +\sinh\theta_{n-1\over 2}  \bbt e_{n-2},\ e_n
\}  \qquad {\rm if\ } n {\rm \ is\ odd,}
$$
}
where $e_1,...,e_n$ is a  space-like $\bbd$-unitary basis for $\bbd^n$ and
 $\theta_1\geq \cdots \geq \theta_{\left[{n\over 2}\right]}\geq 0$.

\pf The case $n=2$ has already been proven.  The 2-form $\o\bigr|_P$ can be put
in canonical form with respect to the positive definite inner product $\langle\cdot,\cdot\rangle$ on $P$.
(This is the same as the canonical form for a skew $n\times n$-matrix under the conjugate action of 
the orthogonal group.) Namely, there exist an orthonormal basisi $\e_1,...,\e_n$ for $P$
and  $\lambda_1\geq \cdots \geq \lambda_{\left[{n\over 2}\right]}\geq 0$ such that:
$$
\o\bigr|_P\ =\ \sum_{j=1}^{\left[{n\over 2}\right]} \lambda_j \e_{2j-1}\wedge \e_{2j}.
\eqno{(A.2)}
$$
Note that span$\{ \e_1,\e_2, \bbt\e_1,\bbt\e_2\}  \cong \bbd^2$
and span$\{ \e_3,\e_4, \bbt\e_3,\bbt\e_4\}  \cong \bbd^2$ are orthogonal
in $\bbd^n$ since $\langle \e_i, \bbt \e_j\rangle = \o(\e_i,\e_j)$.
Thus we can apply induction.\qed
\medskip

Next we keep track of orientations.
Let $\x_0 =   {\overrightarrow\bbr}^n$ denote the space-like $n$-plane
$\bbr^n\ss\bbd^n= \bbr^n\oplus \tau\bbr^n$ equipped with the standard orientation.
Let $\Gsp$ denote the connected component of the Grassmannian 
of oriented space-like $n$-planes in $\bbd^n$  which contains $\x_0$. 
Given $\x\in \Gsp$ the orthonormal basis $\e_1,...,\e_n$ for $P$, described in Proposition A.2,
satisfies $\x=\pm \e_1\wedge\cdots\wedge \e_n$.  If $n$ is odd and the minus sign occurs, we
can replace $e_n=\e_n$ by its negative and have a canonical oriented basis.  If $n$ is even and the 
minus sign occurs, replace $e_n$ and $\theta_{\left[{n\over 2}\right]}$ by their negatives.
This replaces $\e_n =  \cosh \theta_{n\over 2} e_n +\sinh\theta_{n\over 2}  \bbt e_{n-1}$
by its negative, and we have proved the following oriented version of the canonical form.

\def\e{e}
\Prop{A.3 (Canonical Form - Oriented Version)}
 {\sl
Each $\x\in \Gsp$ can be put in the canonical form
$$
\xi\ =\ \e_1 \wedge (\cosh \theta_1 \e_2 +\sinh\theta_1 \bbt \e_1) \wee
\e_{n-1} \wedge ( \cosh \theta_{n\over 2} \e_n +\sinh\theta_{n\over 2}  \bbt \e_{n-1})
 \qquad {\rm if\ } n {\rm \ is\ even, \ and}
$$
$$
\xi\ =\ \e_1\wedge( \cosh \theta_1 \e_2 +\sinh\theta_1 \bbt \e_1) \wee
\e_{n-2} \wedge (\cosh \theta_{n-1\over 2} \e_{n-1} +\sinh\theta_{n-1\over 2}  \bbt \e_{n-2})\wedge  \e_n
  \quad {\rm if\ } n {\rm \ odd,}
$$
where $\e_1,...,\e_n$ is a  space-like $\bbd$-unitary basis for $\bbd^n$ and
 $\theta_1\geq \cdots \geq \theta_{\left[{n\over 2}\right]}$; all $\geq 0$ if $n$ is odd, and all but one
 can be taken $\geq0$ if $n$ is even.  Moreover, $\e_1\wee  \e_n\in \Gsp$.
Consequently, the  map $A\in M_n(\bbd)$ sending the standard space-like unitary  basis 
  $\be_1,...,\be_n$  for $\bbd^n$ to the space-like unitary  basis 
 $\e_1,...,\e_n$  for $\bbd^n$ has  $\det_\bbd A = re^{\tau\theta}$,  $r>0$, that is, 
 $A\in {\rm U}^+_n(\bbd)$.
 }

\pf
The only thing left to prove is that $\e_1\wee \e_n\in \Gsp$. The homothety connecting $\theta$ to 0 connects $\x$ to 
$\e_1\wee \e_n$.  Since $\xi \in \Gsp$, this proves that $\e_1\wee \e_n \in \Gsp$.\qed

\def\e{\epsilon}

 \vskip .3in

\centerline{\headfont Appendix B }  
\medskip

\centerline{\headfont \ Projective and Hermitian Projective varieties over D. }  
\medskip
Much of hermitian projective geometry carries over to the $\bbd$ setting,
with some interesting twists. We begin with the following.

\Def{B.1}  The {\sl $n$-dimensional projective space} $\bbp^n(\bbd)$ over $\bbd$
is the set of all free rank-one submodules of $\bbd^{n+1}$.  Otherwise said,
$\bbp^n(\bbd)$ is the set of all $T$-invariant real 2-dimensional subspaces  $\ell\ss\bbd^{n+1}$
for which the $+1$ and $-1$ eigenspaces $\ell_{\pm}$ are both one-dimensional.

\medskip
Passing to null coordinates    and writing  $\bbd^{n+1} = e\,\bbr^{n+1} \oplus \bar e \, \bbr^{n+1}$
(the $\pm$-eigenspaces of $T$), we see that the decomposition $\ell=\ell_+\oplus \ell_-$ gives
$$
\bbp^n(\bbd)\ =\ \bbp^n(\bbr)\times\bbp^n(\bbr)
 $$
  There is a corresponding decomposition for all the $\bbd$-grassmannians.  In fact we have 
  the following.  For every real projective algebraic variety $V(\bbr)$ there is a corresponding
  variety $V(\bbd)$ over $\bbd$ defined by the 
  same polynomials with arguments now in $\bbd^{n+1}$
  (base change).
Taking null coordinates $z=eu+\bar e v \in \bbd^{n+1} = e\,\bbr^{n+1} \oplus \bar e \, \bbr^{n+1}$,
 a polynomial $p(x)$ with real coefficients satisfies
$
p(z)\ =\ e\, p(u) + \bar e\, p(v)
$
which shows that 
$$
V(\bbd)\ =\ V(\bbr)\times V(\bbr).
\eqno{(B.1)}
$$

Note that in choosing {\sl homogeneous coordinates}  for $\bbp^n(\bbd)$ one must restrict
to elements $z\in \bbd^{n+1}$ which generate a {\sl free} submodule.  This means that
$z= eu+\bar e v$ where $u\neq 0$ and $v\neq 0$. The space $\bbp^n(\bbd)$ is then 
obtained by dividing by the group $\bbd^*$.

When we introduce the split Fubini-Study metric this all becomes more interesting.
Recall the hermitian metric
$(z,\zeta) = \sum_{k=0}^n z_k {\bar \zeta}_k$ with associated quadratic form
$$
(z,z) \ =\  \sum_{k=0}^n u_k {v}_k\ \equiv \  u\bullet v.
$$
In $\bbp^n(\bbd)$ there is the {\sl projective Stiefel variety}:
$$
\eqalign
{
{\rm St} \ &\equiv \   \{[z] \in  \bbp^n(\bbd) : (z,z)=0\}\ \cr
&=\   
  \{[e\,u + \bar e\,v ] \in  \bbp^n(\bbd) \ :\ u\bullet v=0\}  \cr
&\cong \ \{ (\ell_1 , \ell_2) \in \bbp^n(\bbr)\times \bbp^n(\bbr) : \ell_1\perp\ell_2\}
}
$$
which is a hypersurface (dimension $2n-1$) in $\bbp^n(\bbr)\times \bbp^n(\bbr)$.
It is not difficult to check that $(\bbd^{n+1})^*\equiv \bbd^{n+1} - \{(z,z)=0\} =  \bbd^{n+1} - \{u\bullet v=0\}$ has four connected
components corresponding to the four components of the unit sphere $\Sph\equiv \{u\bullet v=1\}$.
They are permuted by the group $\{\pm1\}\times \{\pm\tau\}$. We define the $\bbd$ {\sl hermitian 
projective space} to be the connected open set
$$
\bbp^n_{\rm herm}(\bbd) \ \equiv\ \bbp^n(\bbd) - {\rm St}
$$
and note the mapping
$
(\bbd^{n+1})^* \to  (\bbd^{n+1})^* /\bbd^* =  \bbp^n_{\rm herm}(\bbd)
$
obtained by dividing by the group $\bbd^*$.
Restricting to the unit sphere gives the  {\sl split Hopf mapping}
$$
\Sph\ \arr \  \bbp^n_{\rm herm}(\bbd)
$$
with fibres which are orbits of the  unit circle $S \equiv \{z\in \bbd^+: \bra zz= \pm 1\}$ acting
by scalar multiplication.

We point out that the group ${\rm GL}_{n+1}(\bbd)$ acts transitively
on $\bbp^n(\bbd)$ and the group ${\rm U}_{n+1}(\bbd)$ acts transitively 
on the open subset $\bbp^n_{\rm herm}(\bbd)$. Thus,
$$
\bbp^n(\bbd)\ =\ {\rm GL}_{n+1}(\bbd)/{\rm GL}_{1,n}(\bbd)
\qquad\supset\qquad
\bbp^n_{\rm herm}(\bbd)\ =\ {\rm U}_{n+1}(\bbd)/{\rm U}_{1}(\bbd)\times {\rm U}_{n}(\bbd)
$$
where ${\rm GL}_{1,n}(\bbd)$ is the subgroup preserving the initial line.

Now in hermitian homogeneous coordinates 
 $(\bbd^{n+1})^*$
 consider the real   ${\rm U}_{n+1}(\bbd)$-invariant 2-form
 $$
 \wt\o\ \equiv\  - \tau \partial\dbar \log (z,z) 
 $$
which in null coordinates can be written
$$
 \wt\o\ =\ d_ud_v \log  (u\bullet v)\ =\ {1\over u\bullet v} \sum_jdu_j\w dv_j - 
  {1\over( u\bullet v)^2} (v\bullet du)\w (u\bullet dv).
$$

   Fix a point $(u_0,v_0)\in \bbd^{n+1}$ with $u_0\bullet v_0=1$.  Then $\nu= (u_0,v_0)$,
   translated to $(u_0,v_0)$ is the unit normal to the sphere $\Sph$ at this point, and the 
   {\sl horizontal subspace}   $H_{(u_0,v_0)} \ss T_{(u_0,v_0)}\Sph$ \ \  is the 
   $\bbd$-orthogonal complement of $\nu$. Note that
   $
   \nu \hk \wt\o\ =\ u_0\bullet dv-v_0\bullet du
    - (v_0\bullet u_0)u_0\bullet dv +(u_0\bullet v_0) v_0\bullet du
    \ =\ 0.
   $
It follows that $(T\nu)\hk \wt\o=0$, i.e., that $\wt\o\bigr|_{\Sph}$ is a horizontal, 
${\rm U}_{n+1}(\bbd)$-invariant 2-form. It is straightforward to check that
$\wt\o$ is non-degenerate on $H_{u_0,v_0}$. (It is particularly clear at points $(u_0,u_0)$. Then
use the transitivity of ${\rm U}_{n+1}(\bbd)$.) We have proved the following.

\Lemma {B.2}  {\sl
The form $\wt\o$ descends to a $d$-closed ${\rm U}_{n+1}(\bbd)$-invariant 
symplectic 2-form $\o$ on $\bbp^n_{\rm herm}(\bbd)$.  This form makes 
$\bbp^n_{\rm herm}(\bbd)$ a $\bbd$-\Kahler manifold.
}
\medskip

From Note \G.5 we have the following.
\Cor{B.3} {\sl
Let  $X\ss \bbp^n(\bbd)$  be a $\bbd$-submanifold and 
let $X_{\rm herm} \ss\bbp^n_{\rm herm}(\bbd)$ denote the set of points
in $X$  where $\o\bigr|_X$  is non-degenerate.
Then $X_{\rm herm}$ is a $\bbd$-\Kahler manifold in the induced metric.
}

\medskip

Consider now a $\bbd$-submanifold 
$$
V(\bbd)\ss \bbp^n(\bbd)
\eqno{(B.2)}
$$
 of $\bbd$-dimension $k$
arising from an algebraic submanifold 
$$
V(\bbr) \ss \bbp^n(\bbr)
\eqno{(B.3)}
$$
as above, and recall the canonical product structure (B.1).  
The embedding (B.2) induces a riemannian metric $g$ on $V(\bbr)$ from 
the standard spherical metric $g_0$ on $ \bbp^n(\bbr)$. 
The embedding (B.3)
induces a K\"ahler metric with K\"ahler form $\o$ on the subset $V(\bbd)_{\rm herm}$. There 
is also a nowhere null holomorphic $k$-form 
$$
\Phi\ =\ e\cdot \a + \bar e\cdot \b
$$
where $\a$ and $\b$ denote the volume form of $V(\bbr)$, with metric $g_0$,
pulled back from the first and second factors of (B.1) respectively.  Note that 
$$
{\rm Re}(\Phi) \ =\ \half(\b+\a)  \and
 {\rm Im}(\Phi) \ =\ \half(\b-\a).
 $$
The following pretty geometry is straightforward to verify.
\def\vd{V_{\Delta}}

\Prop{B.4}  {\sl Consider the diagonal embedding $\Delta:V(\bbr) \hookrightarrow V(\bbd)$ given by 
 $\Delta(u)= (u,u)$ with respect to (B.1).  Set $\vd = \Delta (V(\bbr))$. Then:
 \smallskip
 
 (1)\ \ $\vd\ss V(\bbd)_{\rm herm}$,
 
 \smallskip

 (2)\ \  $\vd$ is space-like Lagrangian, in particular 
  \smallskip\centerline{$ \o\bigr|_{\vd} \ =\ 0,$}
 
 \smallskip

 (3)\ \  $\vd$ satisfies the ``special'' condition 
 \smallskip\centerline{$ {\rm Im}(\Phi) \bigr|_{\vd} \ =\ 0,$} 
 
 \smallskip

  (4)\ \ the restriction of the (real part of) the K\"ahler metric to $\vd$ coincides with $ 2 g$,
   \smallskip

 (5)\ \ there is a constant $c=c(k)$ such that: 
  \smallskip\centerline{$  {\rm Re}(\Phi) \bigr|_{\vd} \ =\ c \, {\rm vol}(g).$}
 
 \smallskip
 
 The story is similar for the anti-diagonal embedding $\wt{\Delta}:V(\bbr) \hookrightarrow V(\bbd)$ given by  $\wt{\Delta} (u)=(u,-u)$
except that $V_{\wt{\Delta}}$ is purely time-like instead of space-like, and the induced metric is 
$-2g$.
 }

\medskip

Of course the K\"ahler metric on $V(\bbd)_{\rm herm}$ is not Ricci-flat.  However, it is interesting 
to ask whether the metric and the holomorphic $n$-form $\Phi$ can be altered outside $\vd$ so that
so that $\Phi$ becomes parallel.

  \vskip .3in

\centerline{\headfont Appendix C }  
\medskip

\centerline{\headfont \ Degenerate Projections. }  
\medskip

Many of the results in Special Lagrangian geometry which appeared in [HL$_1$] carry
over to the split case.  One instance of this concerns submanifolds with degenerate
projections. For what follows the reader is referred to  [HL$_1$] for full details and 
proofs.

Decompose $\rn$ as $\rn = \bbr^p\oplus \bbr^{n-p}$ and consider a $p$-dimensional 
submanifold $M$ in $\bbr^p \oplus \tau  \bbr^{n-p}$ along with a smooth function $h$ on 
$M$.  Let $N_{\rm twist}(M)$ denote the normal bundle to $M$ twisted by the antipodal
map on the fibres, defined precisely by (C.2) below.
 This is a Lagrangian submanifold of $\dn$ as is its affine translate by $dh$,
$$
X\ \equiv N_{\rm twist}(M) + dh.
\eqno{(C.1)}
$$
The case of interest is when $M$ is space-like in $\bbr^p \oplus \tau \bbr^{n-p}$.
Then $M\equiv $ graph$(u)$ with $u:\bbr^p \to \bbr^{n-p}$, and we can assume 
$h:\bbr^p \to \bbr$.

Using the decomposition of $\rn$ to decompose the coordinates $z=(x',x''), y=(y',y'')$ the manifold
$X$ is parameterized by sending $x\in \rn$ to 
$$
\left (x', \, x'', \, x''\cdot {\partial u\over \partial x'}(x') +{\partial h\over \partial x'}(x'), \, u(x')       \right).
\eqno{(C.2)}
$$
Therefore, $X$ is the graph (over $\rn$ in $\dn = \rn \oplus \tau \rn$) of $\nabla \vf$, the
gradient of the potential function 
$$
\vf(x) \ \equiv\ x''\cdot u(x') + h(x').
\eqno{(C.1)}
$$
In the case $h\equiv 0$ we have:

\Theorem {C.1} {\sl
If $M$ is a space-like austere submanifold of $\bbr^p\oplus \tau \bbr^{n-p}$, then
$N_{\rm twist}(M)$ is a special Lagrangian submanifold of $\dn$.
}

\medskip

The case $n=3$ and $p=2$ is particularly nice even with $h\neq0$.  Note that
$\vf(x) =x_3u(x_1,x_2) +h(x_1,x_2)$ is the potential function for $X$.

\Theorem {C.2} {\sl
The submanifold $X \equiv N_{\rm twist}(M) + dh$ is special Lagrangian in $\bbd^3$ if and only
if $M^3$ is space-like in $\bbr^2 \oplus \tau \bbr$ with vanishing mean curvature
and $h$ is a harmonic function on $M$.
 }

\pf
If $\vf(x) = x_3u(x_1,x_2) +h(x_1,x_2)$, then 
$\nabla \vf = (x_3 u_1+h_1, x_3 u_2+h_2, u)$, and
$$
\Hess\,\vf \ =\ 
\left(
\matrix
{
x_3 u_{11}+h_{11}  &  x_3 u_{12}+h_{12}  &         u_1  \cr
x_3 u_{21}+h_{21}  & x_3 u_{22}+h_{22}   & u_2  \cr
u_1 & u_2 & 0\cr
}
\right).
$$
Therefore
$$\eqalign
{
\Delta \vf + \det \Hess\,\vf \ =\ &x_3\left[ (1-u_2^2)u_{11} + 2u_1u_2 u_{12} +(1-u_1^2)u_{22}    \right]  \cr
&\quad \quad + \left[  (1-u_2^2)h_{11} + 2u_1u_2 h_{12} +(1-u_1^2)h_{22}    \right]. 
\qquad\qquad\mathqed}
$$

\Remark{C.3} Each such special Lagrangian manifold $X \equiv N_{\rm twist}({\rm graph}(u))$
has a potential $\bar\vf(u_1,u_2,u_3)$ with $X\equiv {\rm graph} \nabla\bar \vf$ in $e\bbr^3\oplus\bar e\bbr^3$.
The function $\bar \vf$ is a solution to:
$$
\bar\vf \ \ {\rm convex\ and\ \ }  \det \Hess \,\bar \vf \ =\ 1.
\eqno{(C.2)}
$$
It would be interesting to examine some explicit examples such as
$y_3 = u(x_1,x_2) = {1\over a} \arccos(ar+c), \ r=\sqrt{x_1^2+x_2^2}$

  \vskip .3in

\centerline{\headfont Appendix D }  
\medskip

\centerline{\headfont \ Singularities and Semi-Riemannian Calibrations. }  
\medskip
It is straightforward to extend the basic set-up of Section \TT\  to manifolds.
Suppose that $X$ is a semi-riemannian manifold of signature $p,q$ with $p,q\geq 1$.
We identify the Grassmannian $G_{\rm space}(p)_x$ of oriented space-like $p$-planes in $T_xX$ 
with a subset of $\L^p T_xX$ by choosing an oriented orthonormal basis for each such $p$-plane
and taking the wedge product. Since $p,q\geq 1$, each $G_{\rm space}(p)_x$  has two connected
components.  A {\bf space-like orientation} on $X$ is a continuous  choice of one the components,
denoted by $G_{\rm space}^+(p)$, and referred to as the Grassmannian of {\bf positive} (or ``future'')
oriented space-like $p$-planes.

\Def{D.1} The {\bf mass ball} $B_x \ss\L^p T_xX$ is the convex hull of  $G_{\rm space}^+(p)_x$.
The  {\bf comass ball} $B_x^* \ss \L^p T_x^*X$ is defined to be the polar:
$$
B_x^*\ \equiv\ \left \{\phi\in \L^p T_x^*X : \phi(\x)\geq 1 \ \ {\rm for\ all\ }\x\in G_{\rm space}^+(p)_x\right\}
$$
(or equivalently for all $\x\in B_x$).  The  {\bf mass cone} $\cc_x$  is the (convex) cone on $B_x$
with vertex at the origin in $\L^p T_xX$, and the  {\bf comass cone} $\cc_x^*$  is the (convex) cone on $B_x^*$
with vertex at the origin in $\L^p T_x^*X$.\medskip

Note that the closure $\overline{\cc_x^*}$ is the polar cone of  $\overline{\cc_x}$,  i.e., $\phi\in 
\overline{\cc_x^*} \iff  \phi(\x)\geq0$   for all $\x\in \overline{\cc_x}$.  Given $\x\in \cc_x$, there exists a 
unique $\l>0$ such that $\l\x \in \partial B_x$. We define
$$
\|\x\|\ \equiv\ {1\over \l} \quad{\rm to\ be\ the\ {\bf mass\ norm\ } of\ }\x
\eqno{(D.1)}
$$

Similarly, the {\bf comass norm} of $\phi\in \cc_x^*$, denoted by $\|\phi\|^*$, is defined
so that $\|\phi\|^*=1 \iff \phi\in\partial B_x^*$ and so that $\|\phi\|^*$ is homogeneous 
of degree 1 for positive scalars.  (Note that in the analogous positive definite case
the mass and comass cones are all of $\L^p T_xX$ and $\L^p T_x^*X$.)
Both $\|\cdot\|$ and $\|\cdot\|^*$ are superadditive on $\cc_x$ and $\cc^*_x$ respectively.

\Def{D.2} A $d$-closed $p$-form $\phi$ on $X$ which takes its values in the comass cone
$\cc$ is a {\bf calibration} if it has comass one at each point, i.e., 
$$
\phi(\x)\ \geq\ 1 \fa \x\in G_{\rm space}^+(p).
\eqno{(D.2)}
$$
A calibration $\phi$ is said to {\bf calibrate} the set $G(\phi) \equiv \{\phi=1\}\cap G_{\rm space}^+(p)$
of $p$-planes in $G_{\rm space}^+(p)$ where $\phi$ attains its minimum.  An oriented submanifold
$M$ of $X$ is a {\bf $\phi$-submanifold} if $\oa M \in G(\phi)$ at each point, or equivalently,
$M$ is positive space-like and $\phi\bigr|_M$ is the unit volume form for $M$.
\medskip

The fundamental theorem for semi-riemannian calibrations (see Theorem \TT.3 for the case
$X=\bbr^{p,q}$) can be stated for submanifolds as follows.  See the generalization allowing 
singularities, Theorem D.10 below, for the proof.

\Prop{D.3} {\sl
Suppose $\phi$ is a calibration of degree $p$ on a  semi-riemannian manifold with a space-like orientation
an of signature
$p,q$ with $p,q\geq1$.  If $M$ is a compact $\phi$-submanifold (possibly with boundary), then
$$
{\rm vol}(M)\ \geq\ {\rm vol}(N)
\eqno{(D.3)}
$$
for all positive space-like compact submanifolds $N$ which are homologous to $M$.  Moreover,
if equality occurs, then $N$ is also a $\phi$-submanifold.
} 

\medskip
The purpose of this appendix is to generalize this result to allow singularities.
\vskip.2in
\centerline{\headfont Currents of Special Type $\cs$.}
\medskip

A current $T$ of dimension $p$ on a manifold $X$ which is {representable by integration}
(see [F]) can be put in {\bf polar form} in each coordinate system $U\ss \rn$
$$
T\ =\ \oa T \|T\|
\eqno{(D.4)}
$$
where $\|T\|$ is a non-negative Radon measure called the {\sl generalized volume measure 
for $T$}, and $\oa T$  is a $\|T\|$-measurable function taking values in the boundary 
of the unit mass ball in $\L^p\rn$, i.e., $|\oa T|=1$ $\|T\|$-a.e..  It is referred to as the {\sl
generalized tangent space} ($\oa T(x)$ may not be simple).  Of course 
$$
T(\vf) \ =\    \int \vf(\oa T)  \|T\|
\eqno{(D.5)}
$$
for each test $p$-form $\vf$ on $U$.

Each riemannian metric on $X$ determines a second notion of mass norm $|\cdot |'$ on 
$\L^pT_xX$ (with $x\in U$).  Defining
$$
 \|T\|' \ \equiv\ | \oa T |' \|T\|\and  { \oa T}' (x) \equiv \ {1\over | \oa T (x)|'} \oa T(x)
\eqno{(D.6)}
$$
yields the polar form
$$
T\ =\ \oa T' \|T\|'
\eqno{(D.7)}
$$
based on the riemannian metric.  Here $\|T\|'$ is a Radon measure on $X$ and
${\oa T}'(x)$ is a $\|T\|'$-measurable section of the bundle $\L^pTX$ with unit mass norm almost everywhere.

Now given a riemannian metric we drop the prime from the notation and use $\oa T$ and $\|T\|$
based on the riemannian metric.

Suppose now that $\cs$ is any cone set (i.e., $s\in \cs_x \Rightarrow ts\in\cs_x \ \forall\, t > 0$)
in the total space of the bundle $\L^pTX$. Given a current $T$ of dimension $p$ which is representable
by integration, if 
$$
 \oa T (x) \in \cs_x \ \ {\rm for} \ \  \|T\|\  {\rm a.a.}\ x 
\eqno{(D.8)}
$$
is true for one riemannian metric,  then it is true for all riemannian metrics.
Consequently, the following concept is independent of the choice of riemannian metric on $X$

\Def{D.4}  Suppose $\cs$ is a cone set in $\L^pT_xX$.  A current $T$ of dimension $p$ which
is representable by integration is said to be of {\bf type $\cs$} if
$$
\oa T(x) \in \cs_x  \ \ {\rm for} \ \  \|T\| \ {\rm a.a.}\ x.
$$

\vskip.2in
\centerline{\headfont The Fundamental Theorem for Semi-Riemannian Calibrations $\cs$.}
\medskip

Suppose that $X$ is a semi-riemannian manifold of signature $p,q\geq1$ which is space-like
oriented and with mass cone $\cc\ss\L^p TX$.

\Def{D.5}    A current $T$ of dimension $p$ which
is representable by integration is said to be {\bf positive space-like} if
$T$ is of type $\cc$.

\Prop{D.6} {\sl
Each positive space-like current $T$ has a {\sl semi-riemannian polar form}
$$
T\ = \ \oa T \|T\|
$$
where $\|T\|$ is a non-negative Radon measure on $X$ and $\oa T$ is a 
$\|T\|$-measurable section of $\L^pTX$ with the semi-riemannian mass norm
$|\oa T(x)| =1$ for $\|T\|$-a.a. points $x$.
}

\pf
Let $T= \oa T' \|T\|'$ denote a local riemannian polar form.  Since $T$ is type $\cc$, 
i.e., ${\oa T}'(x) \in \cc_x$, there is a well defined semi-riemannian mass norm $|{\oa T}'(x) |$. Set
$$
  { \oa T} (x) \equiv \ {1\over | \oa T (x)|'}{ \oa T}'(x)  \and  \|T\| \ \equiv\ | \oa T |' \|T\|'.
\eqno{(D.9)}
$$
Then $T=\oa T\|T\|$ and the semi-riemannian mass norm $|\oa T(x)|$ equals 1 for $\|T\|$-a.a. points $x$.\qed

\Def{D.7} The{\bf volume/mass}  of a positive space-like current $T$  is defined by
$$
{\rm vol}_K(T)\ =\ \int_K \|T\| \qquad{\rm for\ all\ compact\ }K\ss X
$$
using the generalized volume measure given by Proposition D.6.\medskip

\Remark{D.8}  We leave it to the reader to show that this definition of mass for positive
space-like currents $T$ agrees with the definition in [KMW].  Namely, $\int \|T\|$ equals
the infimum of $T(\vf)$ taken over all forms $\vf$ of type $\cc^*$ with comass
$\|\vf\|^*_K \equiv \inf_{x\in K}\|\vf\|^* \geq 1$. (Hint: Use (\TT.3).)

\medskip

Now suppose that $\phi$ is a  calibration of degree $p$ on $X$.

\Def{D.9} A $\phi$-{\bf  subvariety} is a locally rectifiable current $T$ of type $\cc$,
or equivalently of type $\bbr^+\cdot  G^+_{\rm space}(p)$ (the cone on $ G^+_{\rm space}(p)$).
The  $\phi$-{subvariety} is {\bf without boundary} if $dT=0$.

The fundamental theorem of semi-riemannian calibrations can be stated a follows.

\Theorem{D.10} {\sl
Suppose $\phi$ is a calibration of degree $p$ on a space-like oriented, semi-riemannian manifold
of signature $p,q\geq1$.  Suppose that $T$ is a $\phi$-subvariety with compact support.
Then $T$ is homologically volume maximizing.  That is, for any other positive space-like rectifiable current
$S$ with compact support which is homologous to $T$ ($S=T+dR$ with $\supp (R)$ compact), one has
$$
{\rm vol}(T)\ \geq\ {\rm vol}(S)
$$
Moreover, equality occurs if and only if $S$ is also a $\phi$-subvariety.
}

\pf
Since $d\phi=0$, $S=T+dR$, and $\phi$ is a calibration, we have
$$
{\rm vol}(T) \ =\ \int  \|T\| \ =\ \int \phi(\oa T) \|T\| \ =\ T(\phi)
\ =\ S(\phi)\ =\ \int\phi(\oa S) \|S\|\ \geq\ \int  \|S\| \ =\ {\rm vol}(S)
$$
and equality occurs if and only if $\phi(\oa S)=1$ for $\|S\|$-a.a. points.\qed

\Remark{D.11} We leave it to the reader to extend this theorem to currents representable by
integration but not necessarily rectifiable, which are of type $\cc(\phi)$ where
$\cc(\phi)$ is the convex cone on $G(\phi)$.

\Remark{D.12. (Null Planes and the Mass Cone)}  
Let $\overline\cc$ be the closed convex cone on $G^+_{\rm space}(p, \bbr^{p,q})$.
Assume for simplicity that $p\leq q$. Let $\CI(p)$ denote the set of isometries $A:\bbr^p\to\bbr^q$,
i.e., $A^tA=I$, the identity on $\bbr^p$.  Identify $\CI(p)$ with a subset of $\L^p\bbr^{p,q}$ by letting 
graph$(A)$ denote the oriented graph of $A$ and setting $\x_A \equiv (e_1+Ae_1)\wee  (e_p+Ae_p)
\in \L^p\bbr^{p,q}$.  If $p=q=n$ we have the orthogonal group O$(n) \ss\L^n\dn$.

Note that $\bra{e_i+Ae_i}{e_j+Ae_j} = 0$ for all $i,j$ since $A^t A=I$.  That is $\x_A$ is totally null.

\Prop{D.13}{\sl
The closed mass cone $\overline\cc \ss \L^p\bbr^{p,q}$ is the convex cone on $\CI(p)  \ss \L^p\bbr^{p,q}$.
}
\pf  If $A\in \CI(p)$, note that for $0\leq t<1$, $\x_A^t =  (e_1+tAe_1) \wee  (e_p+tAe_p)$
rescaled by $(1-t^2)^{p\over 2}$ is an element of $G^+_{\rm space}(p, \bbr^{p,q})$.
Hence $\CI(p)\ss \overline \cc$.

Conversely, if $\x\in G^+_{\rm space}(p, \bbr^{p,q})$, then pick an oriented orthonormal basis
$e_1,...,e_p$ for $\x$ and an orthonormal set $\bar e_1,...,\bar e_p$ in $\bbr^q$.  Define
$n^\pm_j =\half (e_j\pm \bar e_j)$ so that $e_j= n^+_j + n^-_j$.  Then $\x = e_1\wee e_p
= \sum_{\pm} n_1^\pm \wee n_p^\pm$ and each $2^p  n_1^\pm \wee n_p^\pm \in \CI(p)$.\qed

 \vskip .3in

\centerline{\headfont Appendix E }  
\medskip

\centerline{\headfont \ Singularities in Split SLAG Geometry. }  
\medskip

Recall the definition of type (Definition D.4) and the semi-riemannian polar form (Proposition D.6)
for currents representable by integration.

 \Def{E.1} A   locally rectifiable $n$-current $T$ defined on an open subset of $\dn$

\medskip
\noindent
is  a {\bf split SLAG subvariety} if 
 $$
 \Ta \ \in\ \SLAGt \quad {\rm for \ } \|T\| \ {\rm a.a.\ points,\ }
 $$

\medskip
\noindent
and is  {\bf  positive space-like} if 
 $$
 \Ta \ \in\ \Gsp \quad {\rm for \ } \|T\| \ {\rm a.a.\ points.}
 $$

 The inequality in Theorem \D.4 reads
$$
({\rm Re}\, dz)(\Ta_x) \ \geq\ 1 \qquad{\rm for\ all\  \ }\Ta_x \in \Gsp,
\eqno{(D.2)}
$$
with equality \ $\iff$\  $\Ta_x$ is split special Lagrangian, i.e., $\Ta_x \in \SLAGt$.

Theorem \SS.3 generalizes to include ``singularities''.

\Theorem{E.2}  {\sl
Suppose that $T$ is a split SLAG subvariety with compact support in   $\dn$.
Suppose $S$ is a positive space-like compactly supported rectifiable current on $\dn$
with $dS=dT$.  Then
$$
{\rm vol}(T)\ \geq\ 
{\rm vol}(S)
$$
and  equality holds  if and only if $S$ is also a split  SLAG subvariety.
}
\pf 
This theorem is a special case of the fundamental Theorem D.10. \qed

\medskip
Finally, we note that if the space-like requirement on $\Ta$ is dropped, the currents $T$ are still of interest
(see Appendix F on dimension two).

\Def{E.3}  A   locally rectifiable $n$-current $T$ defined on an open subset of $\dn$
is {\bf unconstrained split SLAG}  (not necessarily space-like split SLAG) if
\medskip

(1) \ \ $T\wedge \o = 0$ \ \ (i.e., $T$ is Lagrangian), and

\medskip

(2) \ \ $T\wedge \psi = 0$ \ \ where $\psi = {\rm Im}\, dz = \half(dv-du)$.

 \vskip .3in

\centerline{\headfont Appendix F }  
\medskip

\centerline{\headfont \ Split SLAG Singularities in Dimension Two. }  
\medskip
This appendix is devoted to describing the structure in dimension $n=2$ of $d$-closed currents
which are unconstrained split SLAG (cf. Definition D.4) and currents that are split SLAG varieties.
We use the corresponding results for complex curves which are collected together as a remark.

\Remark{F.1. (Complex Curves in $\bbc^2$)}  Suppose $V$ is an irreducible complex curve defined near a point
 $z\in \bbc^2$.  Even if $z$ is a singular point of $V$, the curve $V$ has a well-defined tangent line
 (with multiplicity $>1$).  The current $[V]$ of integration over $V$ is well-defined and the singularities
 are isolated.  The Gauss map is an open map which can be understood using a local uniformizing 
 parameter.  A {\bf holomorphic chain} $T$ in an open subset $X$ of $\bbc^2$ is defined to be a locally
 finite sum  of the form $\sum_j n_j [V_j]$ where each $n_j\in \bbz$ and each $V_j$ is a closed 
 one-dimensional subvariety of $X$.  Thus, holomorphic chains and divisors on $X$ are equivalent concepts.
 
 Locally rectifiable 2-currents  $T$  for which $\Ta$ is a field of complex lines, 
  are completely understood because of the 
 following result.
 
 \Theorem{F.2. ([HS], [S])} {\sl
 Suppose that $T$ is a $d$-closed rectifiable 2-current defined in an open subset $X$ of $\bbc^2$.
 If the (unoriented) tangent space $\Ta$ to $T$ is a complex line $\|T\|$-a.e., then $T$ is a holomorphic chain
 in $X$.
 }
\medskip

This result can be viewed as a  very strong regularity result  for ``rectifiable'' complex curves.
Fortunately, it applies directly to split SLAG geometry.

\Theorem{F.3} {\sl
There exists a real orthogonal coordinate change $\bbc^2 \leftrightarrow \bbd^2$, such that holomorphic chains
are transformed to unconstrained split SLAG currents and vice-versa.
}

\pf Let $z'\in \bbc^2$ denote complex coordinates and $z\in \bbd^2$ denote $\bbd$-coordinates
with $z_j'=x_j'+iy_j'$ and $z_j=x_j+\tau y_j$.
Define a coordinate change  $z_1' = x_1+ix_2$ and $z_2'=y_1+iy_2$ from $\bbd^2$ to $\bbc^2$.  Then
$$
dz'\ =\ \o_\bbd + i{\rm Im}\, dz.
\eqno{(F.1)}
$$
Also we have
$$
\o_\bbc\ =\ dx_1 \, dx_2 -dy_1\,dy_2
\and
{\rm Re}\, dz \ =\  dx_1' \, dx_2' -dy_1'\,dy_2'.
$$
Suppose that $T$ is an unconstrained split SLAG current.  Then by definition and (F.1) we have
$T\wedge dz'=0$ for $\|T\|$-a.a. points.  This condition $T\wedge dz'=0$ is equivalent to $\Ta$ being a complex line.
 Hence by the Structure Theorem F.2, $T$ is a holomorphic chain with respect to
the complex structure $J$ on $\bbd^2$ defined by
$$\hskip1in
J\left( { \partial  \over \partial x_1   }  \right) \ =\ { \partial  \over \partial x_2   }
\and
J\left( { \partial  \over \partial y_2   }  \right) \ =\ { \partial  \over \partial y_1   }\hskip 1in\mathqed
$$

Note that the condition 
$$
{\rm Im}\, dz\big|_{{\rm graph}\, A}\ = \ 0\ \qquad{\rm for\ \ } A\in M_2(\bbr)
$$
is that $\tr A=0$.  Hence, $M={\rm graph}\, A$ is unconstrained split SLAG \ $\iff$\ $A$ is symmetric with trace zero.  
Thus 
$$
A\ =\ \left(     
\matrix{
a & b\cr
b&-a\cr
}
\right)
$$
and the graph $y= Ax$ becomes $z_2' = (a+ib)z_1'$ in the complex coordinates 
$z_1', z_2'$   above.

Now $M\equiv {\rm graph}\, A$, with 
$A$ as above, is split SLAG \ $\iff$\ in addition $A^tA =(a^2+b^2)\cdot I < I$, i.e., $$|a+ib|\ <\ 1$$.

\Def{F.4}  
A holomorphic chain $T=\sum_j n_j [V_j]$ in $\bbc^2$ satisfies the {\bf 45$^o$ rule} if 
if each tangent plane is of the form $z_2'=\a z_1'$ with $|\a|<1$ (including tangent planes at singular points).

\Cor{F.5} {\sl
Under the real orthogonal coordinate transformation given in Theorem F.3,
a split SLAG variety $T$ corresponds to a holomorphic chain $T$ which satisfies the 45$^o$ rule.
}

\pf
By the hypothesis on $T$ that $\Ta\in \SLAGt$ at $\|T\|$-a.a. points, we only know that $T$ is a holomorphic
chain whose tangent plane  satisfies the 45$^o$ rule at $\|T\|$-a.a. points.  However, since the Gauss map
of each $V_j$ is open (even in neighborhoods of singular points),
 it follows that then tangent planes satisfy the 45$^o$ rule at every point of the support
of $T$.\qed

\Remark{F.6}
These results can be restated in null coordinates on $\bbd^2$ as follows.  Suppose $T$ is a $d$-closed
locally rectifiable current on an open subset $X$ of $\dn$, with the property that for $\|T\|$-a.a. points
in $X$, the unoriented tangent space of $T$ is the graph of some symmetric matrix $A>0$ with $\det\,A=1$
over $e\bbr^2$ in $\bbd^2 = e\bbr^2\oplus \bar e \bbr^2$.  Then $T$ is a locally finite sum of currents
of integration $\vf_*([\Delta])$ where $\vf:\Delta\to \bbd^2$ maps $\zeta$ to $z(\zeta)$ and 
$$
x_1(\z)+ix_2(\z) \and y_2(\z)+iy_1(\z) \quad{\rm are\ holomorphic}.
\eqno{(F.2)}
$$
Conversely, if $T=\vf_*([\Delta])$ where $\vf$  satisfies (F.2),  and if
$$
\left | {\partial \over \partial \z}\bigl( x_1(\z) +ix_2(\z)\bigr)\right|
\ <\ 
\left | {\partial \over \partial \z}\bigl( y_2(\z) +iy_1(\z)\bigr)\right|,
$$
then, except for isolated singularities,  $\vf$ parameterizes a split SLAG submanifold of $X\ss\bbd^2$.

\vfill\eject



\centerline{\bf References}

\vskip .2in

\noindent
\item{[A$_1$]}  A. D. Alexandrov, {\sl Almost everywhere existence of the second differential of a convex function and properties of convex surfaces connected with it (in Russian)}, 
Lenningrad State Univ. Ann. Math.  {\bf 37}   (1939),   3-35.

 \smallskip

\noindent
\item{[A$_2$]}    A. D. Alexandrov,  {\sl  The Dirichlet problem for the equation Det$\| z_{i,j}\| = \psi(z_1,...,z_n,x_1,...,x_n)$}, I. Vestnik, Leningrad Univ. {\bf 13} No. 1, (1958), 5-24.

\smallskip

\noindent
\item{[AMT]} D. V. Alekseevsky, C. Medori and A. Tomassini, {\sl Homogeneous para-K\'ahler Einstein manifolds}, 
Russian Math. Surveys 64, No. 1, (2009), 1-43. ArXiv:0806.2272.

 \smallskip

\noindent
\item{[B]} R. Bryant, {\sl  Classical, exceptional, and exotic holonomies: a status report}, 
Actes de la Table Ronde de G\'eom\'etrie Diff\'erentielle (Luminy, 1992), 93-165, S\'emin. Congr., 1, Soc. Math. France, Paris, 1996.

 \smallskip

%
%
 
\noindent
 \item{[C$_1$]}   L. A. Caffarelli,    {\sl
Interior a priori estimates for solutions of fully non-linear equations},  
Ann. of  Math.  {\bf 130} (1990),   189-213.

 \smallskip

\noindent
 \item{[C$_2$]}   L. A. Caffarelli,    {\sl
A localization property of viscosity solutions to the Monge-Amp\`ere equation
and their strict  convexity},  
Ann. of  Math.  {\bf 131} (1990),   129-134.

 \smallskip

\noindent
 \item{[C$_3$]}   L. A. Caffarelli,    {\sl
Interior $W^{2,p}$ estimates for solutions of the Monge-Amp\`ere equation},  
Ann. of  Math.  {\bf 131} (1990),    135-150.

 \smallskip

\noindent
 \item{[CNS$_1$]}   L. Caffarelli, L. Nirenberg and J. Spruck,  {\sl
The Dirichlet problem for nonlinear second order elliptic equations. I: Monge-Amp\`ere equation},  
Comm. Pure Appl. Math.  {\bf 37} (1984),   369-402.

 \smallskip

\noindent
 \item{[CNS$_2$]}   L. Caffarelli, L. Nirenberg and J. Spruck,  {\sl
The Dirichlet problem for nonlinear second order elliptic equations, III: 
Functions of the eigenvalues of the Hessian},  Acta Math.
  {\bf 155} (1985),   261-301.

 \smallskip

\noindent
 \item{[CNS$_3$]}   L. Caffarelli, L. Nirenberg and J. Spruck,  {\sl
The Dirichlet problem for the degenerate Monge-Amp\`ere equation},  Rev. Mat. Iberoamericana
  {\bf 2} (1986),   19-27.

 \smallskip

\noindent
 \item{[CNS$_4$]}   L. Caffarelli, L. Nirenberg and J. Spruck,  {\sl  Correction to:
``The Dirichlet problem for nonlinear second order elliptic equations. I: Monge-Amp\`ere equation''},  
Comm. Pure Appl. Math.  {\bf 40} (1987),   659-662.

 \smallskip

\noindent
\item{[CIL]}   M. G. Crandall, H. Ishii and P. L. Lions {\sl
User's guide to viscosity solutions of second order partial differential equations},  
Bull. Amer. Math. Soc. (N. S.) {\bf 27} (1992), 1-67.

 \smallskip

\noindent
\item{[CFG]}   V. Cruceanu, P. Fortuny and P. Gadea,  {\sl
A survey on paracomplex geometry},  
Rocky Mountain J. Math. {\bf 26} (1996), 83-115.

 \smallskip

\item{[D]}  Y. Dong, {\sl On indefinite special Lagrangian submanifolds in indefinite complex Euclidean spaces},
Journal of Geometry and Physics {\bf 59} (2009), 710-726.

 \smallskip

\item {[EST]}  F. Etayo, R. Santamar\'ia amd U. Tr\'ias {\sl The geometry of a bi-lagrangian manifold},
 Diff.  Geom. Appl. {bf 24} (2006). 33-59.

 \smallskip

\noindent
\item{ [Fu]}   L. Fu, {\sl  On the boundaries of Special Lagrangian submanifolds},
Duke Math. J.   {\bf 79}   no. 2 (1995),   405-422.

 \smallskip

 \noindent 
\item {[GM]}   P. Gadea and A. Montesinos Amilibia,  {\sl Spaces of constant para-holomorphic sectional curvature},  Pacific J. Math. {\bf 136} (1989),  85-101.

 \smallskip

 \noindent 
\item {[GGM]}   P. Gadea, J. Grifone  and J. Mu\~noz Masqu\'e,  {\sl Manifolds modelled over free
modules over the double numbers}, Acta. Math. Hungar. {\bf 100} (2003). 187-203.

 \smallskip

 \noindent 
\item {[H]}   F. R. Harvey,  { Spinors and Calibrations},  Perspectives in Math. vol.9, Academic Press,
Boston, 1990.

 \smallskip

 \noindent 
\item {[HL$_1$]}   F. R. Harvey and H. B. Lawson, Jr,  {\sl Calibrated geometries}, Acta Mathematica 
{\bf 148} (1982), 47-157.

 \smallskip

\item {[HL$_{2}$]} F. R. Harvey and H. B. Lawson, Jr., 
{\sl  An introduction to potential theory in calibrated geometry},   Amer. J. Math.    {\bf 131} no. 4 (2009), 893-944.
ArXiv:math.DG/0710.3920.    

\smallskip

\item {[HL$_{3}$]} F. R. Harvey and H. B. Lawson, Jr., {\sl  Duality of positive currents and plurisubharmonic functions in calibrated geometry},    Amer. J. Math.    {\bf 131}  no. 5 (2009), 1211-1240.
  ArXiv:math.DG/0804.1316.
  
\smallskip

\item {[HL$_{4}$]}  F. R. Harvey and H. B. Lawson, Jr., {\sl  Dirichlet duality and the non-linear Dirichlet problem},    Comm. on Pure and Applied Math. {\bf 62} (2009), 396-443.

\smallskip

\item {[HL$_{5}$]} F. R. Harvey and H. B. Lawson, Jr.,  {\sl  Plurisubharmonicity in a general geometric context}, Stony Brook Preprint.  ArXiv 0804.1316.

\smallskip

\item {[HL$_{6}$]} F. R. Harvey and H. B. Lawson, Jr., {\sl  Lagrangian plurisubharmonicity and convexity},  Stony Brook Preprint (2007).

   \smallskip

\item {[HS]}   F.R. Harvey and B. Shiffman,    {\sl  A characterization of
holomorphic chains},    Ann. of Math.,
 {\bf 99}  (1974), 553-587.

\smallskip

\noindent
\item{[Ha$_1$]}   M. Haskins,     
{\sl Constructing special lagrangian cones },    Ph.D. Thesis, Univ. of Texas at Austin, 2000.

   \smallskip

\noindent
\item{[Ha$_2$]}   M. Haskins,     
{\sl Special lagrangian cones },    Amer. J. Math {\bf 126} (2004), 845-871.

   \smallskip

\noindent
\item{[Ha$_3$]}   M. Haskins,     
{\sl  The geometric complexity of Special lagrangian $T^2$-cones },    Invent. Math. {\bf 157} (2004), 11-70.

   \smallskip

\noindent
\item{[HK$_1$]}   M. Haskins and N. Kapouleas,     
{\sl Special lagrangian cones with higher genus links },    Invent. Math. {\bf 167} (2007), 223-294.

   \smallskip

\noindent
\item{[HK$_2$]}   M. Haskins and N. Kapouleas,     
{\sl  Gluing constructions of Special lagrangian cones },    to appear in the ``Handbook of Geometric Analysis''

   \smallskip

\noindent
\item{[Hi$_1$]}   N. Hitchin,     
{\sl The moduli space of special Lagrangian submanifolds},    
Ann. Scuola Norm. Sup. Pisa Cl. Sci. (4) {\bf 25} (1997), 503-515.

   \smallskip

\noindent
\item{[Hi$_2$]}   N. Hitchin,     
{\sl Generalized Calabi-Yau manifolds},    
Quart. J. Math. {\bf 54} (2003), 281-308.

   \smallskip

\noindent
\item{[IZ]}  S. Ivanov and S. Zamkovoy,     
{\sl Parahermitian and paraquaternionic manifolds},    
Diff. Geom. and its Applications {\bf 23} (2005), 205-234.

   \smallskip

\noindent
\item{[J$_1$]}   D. Joyce,  pp. 163-198 in    ``Different  Faces of Geometry'', Int. Math. Ser. (N.Y.), 3, Kluwer/Plenum, New York, 2004.

   \smallskip

\noindent
\item{[J$_2$]}   D. Joyce,     
{\sl U(1)-Invariant  special Lagrangian 3-folds. III. Properties of singular solutions},    Adv. Math. {\bf 192} (2005), 135-182.

   \smallskip

\noindent
\item{[KM]}   Y.-H. Kim and R. J. McCann,
{\sl Continuty, curvature, and the general covariance of optimal transportation},    ArXiv:0712.3077.
To appear in J. Eur. Math. Soc..

   \smallskip

\noindent
\item{[KMW]}   Y.-H. Kim, R. J. McCann and M. Warren,
{\sl Calibrating optimal transportation: a new pseudo-riemannian geometry},    ArXiv:0907.4962

\smallskip

\noindent
\item{[Mc$_1$]}  R. C. McLean, {\sl Deformations of Calibrated Submanifolds},  Duke University Thesis, 1996.

\smallskip

\noindent
\item{[Mc$_2$]}  R. C. McLean, {\sl Deformations of Calibrated Submanifolds},  Comm.  Anal. Geom. {\bf 6}  No. 4 (1998), 705-747.

\smallskip

\noindent
\item{[M$_1$]}  J.  Mealy,     
{Calibrations on semi-riemannian manifolds},    
Ph.D. Thesis, Rice University, 1989.

\smallskip

\noindent
\item{[M$_2$]}  J.  Mealy,     
{Volume maximization in semi-riemannian manifolds},    
Indiana Univ. Math. J. {\bf 40} (1991), 793-814.

\smallskip

\item {[P$_1$]}  A. V. Pogorelov, {\sl   Monge-Amp\`ere equations of elliptic type},
Noordhoff, Groningen, 1964.

\smallskip

\item {[P$_2$]}  A. V. Pogorelov, {\sl   Extrinsic geometry of convex surfaces},
Translations A. M. S., Vol XXXV, 1973.

\smallskip

\item {[RT]} J. B. Rauch and B. A. Taylor, {\sl  The Dirichlet problem for the 
multidimensional Monge-Amp\`ere equation},
Rocky Mountain J. Math {\bf 7}    (1977), 345-364.

\smallskip

\item {[SS]}     L. Sch\"afer and F. Schulte-Hengesbach,    {\sl Nearly pseudo-K\"ahler and nearly
para-K\"ahler six manifolds},    
ArXiv:0912.3271
\smallskip

\smallskip

\item {[S]}      B. Shiffman,    {\sl Complete characterization of holomorphic chains of codimension one},    
Math. Ann. {\bf 274} (1986), 233-256.
\smallskip

\smallskip

\item {[TU]} N. S. Trudinger and J.n I. E. Urbas, {\sl Second derivative estimates for equations of Monge-Amp\`ere type},
Bull. Austral. Math. Soc. {\bf 30}  (1984), 321-334.
\smallskip

\item {[V]}  C. Villani,  {Optimal Transport: Old and New},
Grundlehren der mathematischen Wissenscharften vol. 338, Srpinger, 2009.

\smallskip

\item {[W$_1$]}  M. Warren,  {\sl Calibrations associated to Monge-Ampere equations},
ArXiv:0702291. To appear in Trans. A.M.S..

\smallskip

\item {[W$_2$]}  M. Warren,  {\sl A McLean theorem for the moduli space of Lie solutions to mass transport equations},
ArXiv:1006.1334.

\smallskip

\item {[Wo]}  J. Wood,  {\sl Foliations on 3-manifolds},
 Annals of Math. {\bf 89} 
  (1969),   336-358.

\smallskip

\item {[Y]}  Yu Yuan,  {\sl A priori estimates for solutions of fully nonlinear special lagrangian equations},
 Ann Inst. Henri Poincar\'e  
{\bf 18}  (2001),   261-270.

\vskip.4in

\vbox
{
F. Reese Harvey

Mathematics Department

Rice University 

Houston, TX 77005-1892

\smallskip

e-mail: harvey@rice.edu
}

\vskip .3in

\vbox
{
H. Blaine Lawson, Jr

Mathematics Department

Stony Brook University

Stony Brook, NY11794-3651

\smallskip

e-mail: blaine@math.sunysb.edu
}

\end